\documentclass[11pt]{amsart}
\usepackage{amsmath,bm}
\usepackage{amssymb}
\usepackage{tabularx}
\usepackage{enumerate}
\usepackage{graphicx}
\newtheorem{thm}{Theorem}[section]
\newtheorem{lem}[thm]{Lemma}
\newtheorem{cor}[thm]{Corollary}
\newtheorem{prop}[thm]{Proposition}
\newtheorem{rmk}[thm]{Remark}
\newcommand{\bpf}{{\bf Proof:\ \ }}
\newcommand{\epf}{\mbox{}\hfill $\Box$}

\numberwithin{equation}{section}
\newcommand{\ol}{\overline}

\newcommand{\be}{\begin{equation}}
\newcommand{\ee}{\end{equation}}
\newcommand\bes{\begin{eqnarray}} \newcommand\ees{\end{eqnarray}}
\newcommand{\bess}{\begin{eqnarray*}}
\newcommand{\eess}{\end{eqnarray*}}

\newcommand{\R}{\mathbb{R}}

\begin{document}

\thispagestyle{empty}

\title[The diffusive competition model with a free boundary]{The diffusive competition model with
a free boundary: Invasion of a superior or inferior competitor$^*$}
\thanks{$^*$ This work was supported by the Australian Research Council and also by NSFC 11071209 of China.}
\date{\today}
\author[Y. Du and Z.G. Lin]{Yihong Du$^\dag$ and Zhigui Lin$^\ddag$}
\thanks{$^\dag$ School of Science and
Technology, University of New England, Armidale, NSW 2351,
Australia}
\thanks
{$^\ddag$ School of Mathematical Science, Yangzhou University,
Yangzhou 225002, China. }
\thanks{{\bf Emails:} {\sf ydu@turing.une.edu.au (Y. Du),
zglin68@hotmail.com (Z. Lin)}}

\begin{abstract} In this paper we consider the diffusive competition model consisting of an invasive species with density $u$ and a native species
with density $v$, in a radially symmetric setting with free boundary. We assume that $v$ undergoes diffusion and growth in $\R^N$, and $u$ exists initially
in a ball $\{r<h(0)\}$, but invades into the environment with spreading front $\{r=h(t)\}$, with $h(t)$ evolving according to the free boundary condition
$h'(t)=-\mu u_r(t, h(t))$, where $\mu>0$ is a given constant and $u(t,h(t))=0$. Thus the population range of $u$ is the expanding ball $\{r<h(t)\}$,
while that for $v$ is $\R^N$. 
In the case that $u$ is a superior competitor (determined by the reaction terms),
we show that a spreading-vanishing dichotomy holds, namely, as $t\to\infty$, either $h(t)\to\infty$ and $(u,v)\to (u^*,0)$, or $\lim_{t\to\infty} h(t)<\infty$ and $(u,v)\to (0,v^*)$, where $(u^*,0)$ and $(0, v^*)$ are the semitrivial steady-states of the system. Moreover, when spreading of $u$ happens, some rough estimates of the spreading speed are also given. When $u$ is an inferior competitor, we show that $(u,v)\to (0,v^*)$ as $t\to\infty$.
\end{abstract}

 \maketitle

\section{Introduction}

In this paper we study the behavior of the  solution
$(u(t, r), v(t, r), h(t))$ to the following reaction-diffusion problem with radial symmetry,
\begin{align}
\left\{
\begin{array}{lll}
u_{t}-d_1 \Delta u=u(a_1-b_1 u-c_1 v),\; & t>0, \ 0\leq r<h(t),  \\
v_{t}-d_2 \Delta v=v(a_2-b_2 u-c_2 v),\; & t>0, \ 0\leq r<\infty,  \\
u_r(t,0)=v_r(t,0)=0,\; u(t, r)=0,\quad & t>0,\ h(t)\leq r<\infty,\\
h'(t)=-\mu u_r(t,h(t)),\quad &t>0,\\
h(0)=h_0,  u(0, r)=u_{0}(r),\;  &0\leq r\leq h_0,\\
 v(0, r)=v_0(r), &0\leq r<\infty,\\
\end{array} \right.
\label{f1}
\end{align}
where $\Delta u=u_{rr}+\frac{N-1}{r}u_r,\ r=h(t)$ is the moving
boundary to be determined, $h_0$, $\mu$, $d_i$,
 $a_i$, $b_i$ and $c_i (i=1,2)$ are given positive constants, and the initial functions
$u_0$  and  $v_0$ satisfy
\begin{align}
\left\{
\begin{array}{ll}
u_0\in C^{2}([0, h_0]), \ u_0'(0)=u_0(h_0)=0\ \ \textrm{and} \ u_0>0\ \textrm{in}\ [0, h_0),\\
v_0\in C^{2}([0, +\infty))\cap L^\infty(0, +\infty), \ v_0'(0)=0\ \ \textrm{and} \ v_0\geq 0\
\textrm{in}\ [0, +\infty).
\end{array} \right.
\label{Ae}
\end{align}
Ecologically, this problem describes the dynamical process of a new competitor invading into the habitat of a native species.
 The first
species ($u$), which exists initially in the ball $\{r<h_0\}$, disperses through random diffusion over an expanding ball
 $\{r<h(t)\}$, whose boundary $\{r=h(t)\}$ is  the invading front, and evolves according to the free boundary condition $h'(t)=-\mu u_r(t, h(t))$, where $\mu$ is a given positive constant.
  The second species ($v$) is
native, which undergoes diffusion and growth in the entire available habitat (assumed to be $\R^N$ here).
The constants $d_1$ and $d_2$ are the diffusion rates of $u$ and $v$,
respectively, $a_1$ and $a_2$ are the intrinsic growth rates, $b_1$ and $c_2$ are the intraspecific and
   $c_1$ and $b_2$ the interspecific competition rates.

In the absence of a native species, namely $v\equiv 0$, the system reduces to the following diffusive logistic problem,
\begin{eqnarray}
\left\{
\begin{array}{lll}
u_{t}-d_1 \Delta=u(a_1-b_1 u),\; & t>0, \ 0\leq r<h(t),  \\
u_r(t,0)=0,\; u(t, h(t))=0,\quad & t>0, \\
h'(t)=-\mu u_r(t,h(t)),\quad &t>0,\\
h(0)=h_0,  u(0, r)=u_{0}(r),\; &0\leq r\leq h_0,
\end{array} \right.
\label{f2}
\end{eqnarray}
which has been treated in \cite{DG}, extending the one dimensional case  first studied in \cite{DL}.
The behavior of \eqref{f2} is characterized by a spreading-vanishing dichotomy, namely, as $t\to\infty$, one of the following alternatives occurs:
\begin{itemize}
\item {Spreading:} $h(t)\to \infty$ and $u(t,r)\to a_1/b_1$, or

\item {Vanishing:} $h(t)\to h_\infty<\infty$ and $u(t,r)\to 0$.
\end{itemize}
Moreover, when spreading occurs, it is shown that
$h(t)/t\to k_0\in (0, 2\sqrt{a_1d_1})$ as $t\to\infty$, and $k_0$ is called the asymptotic spreading speed of $u$. Further discussions of $k_0$ and
a deduction of the free boundary condition based on ecological assumptions can be found in \cite{BDK}.

In this paper, we will examine the case that $u$ invades into an environment where a native competitor already exists. This is a much more
complicated situation, and we will only consider \eqref{f1} under certain restrictions on the parameters, to be specified below.

Problem \eqref{f1} is a variation of  the diffusive Lotka-Volterra  competition model, which  is often considered over a bounded spatial domain
with suitable boundary conditions or considered over the entire space $\R^N$ (\cite{CC, P}).
For example, the dynamical behavior of the following bounded domain problem
\bes\label{us21}
 \left\{\begin{array}{ll} u_{t}-d_1\Delta
u=u(a_1-b_{1}u-c_{1}v),& (t,x)\in (0, \infty)\times\Omega,\\
v_{t}-d_2\Delta v=v (a_2-b_{2}u-c_{2}v),& (t,x)\in (0,\infty)\times\Omega,\\
\frac{\partial u}{\partial \eta}=\frac{\partial v}{\partial \eta}=0,\ & (t,x)\in (0,\infty)\times\partial\Omega,\\
u(0,x)=u_0(x)> 0,\ v(0,x)=v_0(x)> 0,& x\in \Omega
\end{array}\right.
\ees
is well known, where $\Omega$ is a bounded smooth domain of ${\mathbb R}^N$ with $N\geq 1$, $\eta$ is the outward unit normal vector on $\partial \Omega$.
This model describes the situation that two competitors evolve in a closed habitat $\Omega$, with no flux  across the boundary $\partial\Omega$.
Therefore their competitive strengths are completely determined by the coefficients $(a_i,b_i,c_i,d_i)$ in the system, $i=1,2$.

Problem (\ref{us21}) admits the trivial steady state $R_0=(0,0)$ and semi-trivial steady-states $R_1=(a_1/b_{1},\ 0)$ and $R_2=(0,\ a_2/c_{2})$.
Moreover, if $b_1/b_2>a_1/a_2>c_1/c_2$ or $b_1/b_2<a_1/a_2<c_1/c_2$, the problem has a unique constant positive steady-state
\begin{align*}
R^*=\left(\frac{a_{1}c_2-a_{2}c_1}{b_{1}c_{2}-
b_{2}c_{1}},\frac{a_{2}b_1-a_{1}b_2}{b_{1}c_{2}-b_{2}c_{1}}\right).
\end{align*}
These are all the nonnegative constant steady-states. There may also exist non-constant positive
 steady-states, but they are all linearly unstable when $\Omega$ is convex (\cite{KW}).
For the constant equilibria, their roles are summarized below (see, for example, \cite{P} page 666):
\begin{enumerate}[(1)]
\item $R_0$ is always unstable;
\item  when $b_1/b_2>a_1/a_2>c_1/c_2$, $R^*$ is globally asymptotically stable;
\item when $a_1/a_2>\max \{b_1/b_2, c_1/c_2\}$, $R_1$ is globally asymptotically stable;
\item when $a_1/a_2<\min \{c_1/c_2, b_1/b_2\}$, $R_2$ is globally asymptotically stable;
\item when $b_1/b_2<a_1/a_2<c_1/c_2$, $R_1$ and $R_2$ are locally asymptotically stable, and $R^*$ is unstable.
\end{enumerate}

In case (2), the competitors co-exist in the long run, and it is often referred to as the weak competition case, where no competitor wins on loses in the competition.
In case (3), the competitor $u$ wipes $v$ out in the long run and wins the competition; so we will call $u$ the superior competitor
and $v$ the inferior competitor. Analogously $u$ is the inferior competitor and $v$ is the superior competitor in case (4).
Case (5) is the strong competition case, and the long-time dynamics of \eqref{us21} is usually complicated and difficult to determine.

We will only consider cases (3) and (4) for \eqref{f1}. We will show that in case (3), similar to \eqref{f2},  a spreading-vanishing dichotomy holds
for \eqref{f1}, namely
as $t\to\infty$, either $h(t)\to\infty$ and $(u,v)\to R_1$ (spreading for $u$), or $h(t)\to h_\infty<\infty$ and $(u,v)\to R_2$ (vanishing for $u$).
Clearly this is strikingly different to the long-time behavior of \eqref{us21}. However, in case (4), we show that as $t\to\infty$,
 $(u,v)\to R_2$, so the dynamical behavior is similar to that of \eqref{us21} in case (4).

For the entire space problem
\bes\label{entire}
 \left\{\begin{array}{ll} u_{t}-d_1\Delta
u=u(a_1-b_{1}u-c_{1}v),& (t,x)\in (0, \infty)\times\R^N,\\
v_{t}-d_2\Delta v=v (a_2-b_{2}u-c_{2}v),& (t,x)\in (0,\infty)\times\R^N,\\
\end{array}\right.
\ees
extensive work has been done concerning the existence of traveling wave solutions in space dimension $N=1$.
For example, for case (3), it is shown in \cite{KO} that there exists $c_*>0$ such that for each $c\geq c_*$,
\eqref{entire} with $N=1$ has a solution of the form
\begin{equation}
\label{tw}
(u(x,t), v(x,t))=(U(x-ct), V(x-ct))
\end{equation}
satisfying
\begin{equation}
\label{twasym}
U'<0, V'>0, \; (U(-\infty), V(-\infty))=R_1,\; (U(+\infty), V(+\infty))=R_2;
\end{equation}
there is no such solution when $c<c_*$. The general long-time behavior of the Cauchy problem of \eqref{entire}, however, is still
poorly understood (see Remark \ref{rmk-u-inf} below for a partial result).

We end the introduction by mentioning some related research.
In \cite{Lin}, a predator-prey model in one space dimension was considered, where the available habitat is assumed to be a bounded interval
$[0,l]$, and no-flux boundary conditions are assumed for both species, except for the predator at $x=l$. It is assumed that
the predator satisfies  a free boundary condition as in \eqref{f1}, before the free boundary $x=h(t)$ reaches  $x=l$, and a no-flux boundary condition at $x=l$
is
satisfied by the predator  after the free boundary has reached $x=l$. It is shown in \cite{Lin} that the free boundary always reaches $l$ in finite time, and hence the long-time dynamical behavior
of the free boundary problem is the same as the fixed boundary problem.
After the first version of this paper was completed, we have learned several more closely related research. In \cite{GW}, the week competition case
was considered in one space dimension, but in their model, both species share the same free boundary.
Such a free boundary setting was also used in \cite{mW} for the Lotka-Volterra predator-prey system in one space dimension. In \cite{WZ},
the Lotka-Volterra predator-prey model was considered in one space dimension, where similar to \eqref{f1}, one species (the predator) is subject to free boundary conditions,
and the other is considered over the entire $\R^1$.

The rest of this paper is organized in the following way. In section 2, we prove some general existence and uniqueness results, which implies in particular that \eqref{f1} has a unique solution defined for all $t>0$. Moreover, some rough a priori estimates are given, as well as a rather general comparison result. These results are useful here and possibly elsewhere.
In section 3, we investigate the case that $u$ is an inferior competitor, namely the coefficients fall into case (4). Sections 4 and 5 are devoted to the case that $u$ is a superior competitor. A spreading-vanishing dichotomy is established in section 4, and a sharp criterion to distinguish the dichotomy is also given there. In section 5, some rough estimates for the spreading speed is given for the case that spreading of $u$ happens.

\section{Preliminary Results}

In this section,  we first
prove a local existence and uniqueness result for a general free boundary problem, and then we obtain
global existence results, which imply that the solution to (\ref{f1})
 exists for all time $t\in (0,\infty)$. Lastly, we obtain some  comparison results, which will be used in the other sections.

 Consider the following general free boundary problem:
\begin{align}
\left\{
\begin{array}{lll}
u_{t}-d_1 \Delta u=f(u, v),\; & t>0, \ 0\leq r<h(t),  \\
v_{t}-d_2 \Delta v=g(u, v),\; & t>0, \ 0\leq r<\infty,  \\
u_r(t,0)=v_r(t,0)=0,\; u(t, r)=0,\quad & t>0,\ h(t)\leq r<\infty,\\
h'(t)=-\mu u_r(t,h(t)),\quad &t>0,\\
h(0)=h_0,  u(0, r)=u_{0}(r),\;  &0\leq r\leq h_0,\\
 v(0, r)=v_0(r), &0\leq r<\infty,\\
\end{array} \right.
\label{1f1}
\end{align}
where $f(0,v)=g(u,0)=0$ for any $u, v\in \mathbb{R}$, and $u_0$, $v_0$ are as in \eqref{f1}.
\begin{thm}\label{localex} Assume that $f$ and $g$ are locally Lipschitz continuous in ${\mathbb R}^2_+$.
For any given $(u_0, v_0)$ satisfying \eqref{Ae} and any $\alpha \in (0,1)$,
there is a $T>0$ such that problem \eqref{1f1} admits a unique bounded
solution $$(u, v, h)\in C^{(1+\alpha)/2,
1+\alpha}({D}_{T})\times C^{(1+\alpha)/2,
1+\alpha}({D}^\infty_{T})\times C^{1+\alpha/2}([0,T]);$$ moreover,
\begin{eqnarray} \|u\|_{C^{
(1+\alpha)/2,
1+\alpha}({D}_{T})}+\|v\|_{C^{
(1+\alpha)/2,
1+\alpha}({D}^\infty_{T})}+\|h\|_{C^{1+\alpha/2}([0,T])}\leq
C,\label{Ba}
\end{eqnarray}
where  $D_{T}=\{(t,r)\in \mathbb R^2: t\in [0,T], r\in [0, h(t)]\}$,
$D^\infty_{T}=\{(t,r)\in \mathbb R^2: t\in [0,T], r\in [0, +\infty)\}$,
$C$ and $T$ only depend on $h_0$, $\alpha$, $\|u_0\|_{C^{2}([0,
h_0])}$, $\|v_0\|_{C^{2}([0,
\infty))}$ and the local Lipschitz coefficients of $f, g$.
\end{thm}
\bpf  The proof is similar to that in \cite{DL} and \cite{DG} for the scalar problem, with some modifications.
We sketch the details here for completeness.
First we straighten the free boundary as in \cite{CF1}. Let
$\zeta (s)$ be a function in $C^3[0, \infty)$ satisfying
\begin{align*}
\zeta (s)=1\ \  \textrm{if} \,\, |s-h_0|<\frac{h_0}8,\; \zeta (s)=0\ \
 \textrm{if} \,\, |s-h_0|>\frac {h_0}2,\quad |\zeta
'(s)|<\frac 5{h_0} \mbox{ for all } s.
\end{align*}
Consider the transformation
\begin{align*}
(t, y)\rightarrow (t, x), \textrm{where}\,\, x=y+\zeta (|y|)(h(t)-h_0)\frac{y}{|y|},
\quad y\in R^N,
\end{align*}
which leads to the transformation
\begin{align*}
(t, s)\rightarrow (t, r), \textrm{with}\,\, r=s+\zeta (s)(h(t)-h_0),
\quad 0\leq s<\infty.
\end{align*}
As long as
\begin{align*}
|h(t)-h_0|\leq \frac {h_0}8,
\end{align*}
the above transformation $x \to y$ is a diffeomorphism from ${\mathbb R}^N$
onto $\mathbb R^N$ and the transformation $s\to r$ is also a diffeomorphism from $[0, +\infty)$
onto $[0, +\infty)$. Moreover, it changes the unknown free boundary $|x|=h(t)$
to the fixed sphere $|y|=h_0$. Now, direct calculations show that
\begin{align*}
\displaystyle\frac {\partial s}{\partial r}=\frac
1{1+\zeta'(s)(h(t)-h_0)}&\equiv
\sqrt{A(h(t), s)},\\
\displaystyle \frac {\partial^2 s}{\partial r^2}=-\frac
{\zeta''(s)(h(t)-h_0)}{[1+\zeta'(s)(h(t)-h_0)]^3}&\equiv
B(h(t),s),\\
\displaystyle -\frac 1{h'(t)}\frac {\partial s}{\partial t}=\frac
{\zeta(s)}{1+\zeta'(s)(h(t)-h_0)}&\equiv C(h(t), s).
\end{align*}
Let us also denote
\[
\displaystyle \frac{(N-1)\sqrt{A}}{s+\zeta(s)(h(t)-h_0)}\equiv D(h(t), s).
\]

If we set
\begin{align*}
&u(t, r)=u(t, s+\zeta (s)(h(t)-h_0))=w(t, s),\\
&v(t, r)=v(t, s+\zeta (s)(h(t)-h_0))=z(t, s),
\end{align*}
then the free boundary problem \eqref{1f1} becomes
\begin{eqnarray}
\left\{
\begin{array}{lll}
w_{t}-Ad_1w_{ss}-(Bd_1+h'C+Dd_1)w_s=f(w,z),&t>0,\; 0\leq s<h_0, \\
z_{t}-Ad_2z_{ss}-(Bd_2+h'C+Dd_2)z_s=g(w,z),&t>0,\; 0\leq s<\infty, \\
w_s(t,0)=z_s(t, 0)=w(t, r)=0, &t>0,\; h_0\leq s<\infty,\\
h'(t)=-\mu w_s(t, h_0), &t>0, \\
h(0)=h_0,\quad w(0, s)=w_0(s):=u_{0}(s), &0\leq s\leq h_0,\\
z(0, s)=z_0(s):=v_0(s),&0\leq s<\infty,
\end{array} \right.
\label{Bb}
\end{eqnarray}
where $A=A(h(t),s)$, $B=B(h(t),s)$, $C=C(h(t),s)$ and $D=D(h(t),s)$.

We denote $h^*=-\mu u_0'(h_0)$, $w_0(s)=0$ for $s>h_0$ and for $0<T\leq\frac {h_0}{8(1+h^*)}$, set
\begin{align*}
&\Delta_{T}=[0,T]\times [0, h_0],\quad \Delta^\infty_{T}=[0,T]\times [0, +\infty),\\
&H_T=\Big\{h\in C^1[0,T]: \, h(0)=h_0,  \ h'(0)=h^*, \, \|h'-h^*\|_{C([0, T])}\leq 1\Big\},\\
&W_T=\Big\{w\in C(\Delta^\infty_{T}): \,w(t,s)\equiv 0 \,\ \textrm{for}\ s\geq h_0, 0\leq t\leq T,\\
&\ \ \ \ \ \ \ \ \ \ \ \ \ w(0,s)=w_0(s)\, \ \textrm{for}\, \ 0\leq s\leq h_0, \
\|w-w_0\|
_{L^\infty(\Delta^\infty_{T})}=\|w-w_0\|
_{C(\Delta_{T})}\leqslant 1\Big\}, \\
&Z_T=\Big\{z\in C(\Delta^\infty_{T}): \, z(0,s)=z_0(s),\,  \|z-z_0\|
_{L^\infty(\Delta^\infty_{T})}\leqslant 1\Big\}.
\end{align*}

It is not difficult to
see that $\Gamma_T:=W_T\times Z_T\times H_T$ is a complete metric space with
the metric
\begin{align*}
\mathcal{D}((w_1, z_1, h_1), (w_2, z_2, h_2))=\|w_1-w_2\|_{C(\Delta_{T})}+\|z_1-z_2\|_{L^\infty (\Delta^\infty_{T})}+
\|h'_1-h'_2\|_{C([0,T])}.
\end{align*}
Let us observe that for $h_1, h_2\in
H_T$, due to $h_1(0)=h_2(0)=h_0$, we have
\begin{equation}
\label{Bc}
 \|h_1-h_2\|_{C([0,T])}\leq  T\|h'_1-h'_2\|_{C([0, T])}.
\end{equation}

Next, we shall prove the existence and uniqueness result by using
the contraction mapping theorem. Since $f$ and $g$ are locally Lipschitz continuous,
there exists an $L^*$ depending on $\|u_0\|
_{C([0, h_0])}$ and $\|v_0\|_{L^\infty ([0, +\infty))}$ such that
\begin{align*}
&|f(w,z)|=|f(w,z)-f(0,z)|\leq L^* |w|\leq L^*(\|u_0\|_{C[0, h_0]}+1), \ (w, z)\in W_T\times Z_T, \\
&|g(w,z)|=|g(w,z)-g(w,0)|\leq L^* |z|\leq L^*(\|v_0\|_{L^\infty[0, \infty)}+1), \ (w, z)\in W_T\times Z_T.
\end{align*}
By  standard $L^p$ theory and the Sobolev imbedding theorem
\cite{LSU},  for any $(w, z, h)\in\Gamma_T$, the
following initial boundary value problem
\begin{align}\label{Bd}
\left\{
\begin{array}{lll}
\tilde w_{t}-Ad_1 \tilde w_{ss}-(Bd_1+h'C+Dd_1)\tilde w_s
=f(w,z),&t>0,\; 0\leq s<h_0, \\
\tilde z_{t}-Ad_2\tilde z_{ss}-(Bd_2+h'C+Dd_2)\tilde z_s=g(w,z),&t>0,\; 0\leq s<\infty, \\
\tilde w_s(t,0)=\tilde z_s(t, 0)=\tilde w(t, r)=0, &t>0,\; h_0\leq r<\infty,\\
\tilde w(0, s)=w_0(s):=u_{0}(s), &0\leq s\leq h_0,\\
\tilde z(0, s)=z_0(s):=v_0(s),&0\leq s<\infty
\end{array} \right.
\end{align}
admits a unique bounded solution $(\tilde{w}, \tilde z)\in C^{(1+\alpha)/2,1+\alpha
}(\Delta_{T})\times C^{(1+\alpha)/2,1+\alpha
}(\Delta^\infty_{T})$ and
\begin{align}
\|\tilde{w}\|_{C^{(1+\alpha)/2,1+\alpha}(\Delta_{T})}\leqslant C_1,\label{Be}
\end{align}
\begin{align}
\|\tilde{z}\|_{C^{(1+\alpha)/2,1+\alpha}(\Delta^\infty_{T})}\leqslant C_1,\label{Be2}
\end{align}
where $C_1$ is a constant depending on
$\alpha, h_0, L^*$,  $\|u_0\|_{C^{2}[0, h_0]}$ and $\|v_0\|_{C^{2}[0, +\infty)}.$
The estimate (\ref{Be2}) comes from the interior estimate. For any $m\geq 0$, by classical parabolic regularity theory \cite{LSU},
one then have the estimate
\begin{align*}
\|\tilde{z}\|_{C^{(1+\alpha)/2,1+\alpha}([0,T]\times(m,m+1)
)}\leq C(\alpha)\|\tilde{z}\|_{W^{2, 1, p}([0,T]\times[m,m+1])}\leqslant C(\alpha, p, L^*, T)
\end{align*}
for some large $p>1$ and $\alpha\in (0, 1)$.

Now, we define $\tilde{h}(t)(>0)$ by the fourth equation in
(\ref{Bb}):
\begin{align}
\label{Bf}
\tilde{h}(t)=h_0-\mu\int^t_0\tilde{w}_s(\tau,h_0)\textrm{d}\tau,
\end{align}
which infers $\tilde{h}'(t)=-\mu \tilde{w}_s(t,h_0)$,
$\tilde{h}(0)=h_0$ and $\tilde{h}'(0)=-\mu
u'_0(h_0)=h^*$. Hence $\tilde{h}'\in C^{\alpha/2}([0,T])$ with
\begin{align}
\label{Bg}
\|\tilde{h}'\|_{C^{(1+\alpha)/2}([0,T])}\leq C_2:=\mu C_1.
\end{align}

Now we define the map
\begin{align*}
\mathcal{F}:\ \Gamma_{T}\longrightarrow C(\Delta^\infty_{T})\times C(\Delta^\infty_{T})\times
C^1[0,T]
\end{align*}
by
$\mathcal{F}(w(t, s), z(t, s); h(t)) = (\tilde w(t, s), \tilde z(t, s); \tilde{h}(t))$.
It's easy to see that $(w(t, s), z(t, s); h(t))\in \Gamma_T$ is a fixed point of
$\mathcal{F}$ if and only if it solves (\ref{Bb}).

The estimates in (\ref{Be}), (\ref{Be2}) and (\ref{Bg}) yield
\begin{align*}
&\|\tilde h'-h^*\|_{C([0,T])}\leq  \|\tilde h'\|_{C^{\alpha /2}([0,T])}T^{\alpha /2}\leq \mu C_1T^{\alpha /2},\\
&\|\tilde w-w_0\|_{C(\Delta _{T})}\leq \|\tilde w-w_{0}\|_{C^{(1+\alpha) /2, 0}(\Delta _{T})}T^{(1+\alpha) /2}\leq
C_1T^{(1+\alpha) /2},\\
&\|\tilde z-z_0\|_{L^\infty(\Delta^\infty_{T})}\leq \|\tilde z-z_{0}\|_{C^{(1+\alpha) /2, 0}(\Delta^\infty_{T})}T^{(1+\alpha) /2}\leq
C_1T^{(1+\alpha) /2}.
\end{align*}
 Therefore if we take $T\leq \min
\{(\mu C_1)^{-2/\alpha},\, C_1^{-2/(1+\alpha)}\}$, then
$\mathcal{F}$ maps $\Gamma_{T}$ into itself.

Now we prove that for $T>0$ sufficiently small, $\mathcal {F}$ is a contraction mapping on
$\Gamma_{T}$. Indeed, let $(w_i, z_i,
h_i)\in \Gamma_{T}$ ($i=1,2$) and denote $(\tilde w_i, \tilde z_i, \tilde
h_i)=\mathcal{F}(w_i, z_i, h_i)$. Then it follows from (\ref{Be}), (\ref{Be2}) and
(\ref{Bg}) that
\begin{align*}
\|\tilde w_i\|_{C^{(1+\alpha)/2, 1+\alpha}(\Delta _{T})}\leq C_1,\, \|\tilde z_i\|_{C^{(1+\alpha)/2, 1+\alpha}(\Delta^\infty_{T})} \leq
C_1,\, \|\tilde h'_i\|_{C^{\alpha/2}([0,T])}\leq C_2.
\end{align*}
Setting $W=\tilde w_1-\tilde w_2$, we find that $W(t, s)$
satisfies
\begin{align*}
& W_t-A(h_2,s)d_1W_{ss}-(B(h_2,s)d_1+h'_2C(h_2,s)+D(h_2,s)d_1)W_s\\
& \;\;\; =[A(h_1,s)-A(h_2,s)]d_1\tilde w_{1,ss}+[B(h_1,s)-B(h_2,s)+D(h_1,s)-D(h_2,s)]d_1 \tilde w_{1,s}\\
& \;\;\;\; +[h'_1C(h_1,s)-h'_2C(h_2,s)]\tilde w_{1,s}+f(w_1, z_1)-f(w_2, z_2),\  t>0, \ 0<s<h_0,  \\
&  \frac {\partial W}{\partial s}(t, 0)=0,\quad W(t, h_0)=0, \quad t>0,\\
&  W(0, s)=0,\quad 0\leq s\leq h_0.
\end{align*}

Using the $L^p$ estimates for parabolic equations and Sobolev's
imbedding theorem, we obtain
\begin{align} \label{e6}
\begin{array}{l}
\|\tilde w_1-\tilde w_2\|_{C^{(1+\alpha)/2,1+\alpha}(\Delta_{T})}  \\
\ \ \ \ \ \ \ \ \ \leq C_3(\|w_{1}-w_{2}\|_{C(\Delta_T)}+\|z_1-z_2\|_{L^\infty(\Delta^\infty_{T})}+\|h_1-h_2\|_{C^1([0,T])}),
\end{array}
\end{align}
where $C_3$ depends on $C_1, C_2$, the local Lipschitz coefficients of $f, g$ and the functions $A, B$ and $C$
in the definition of the transformation $(t,s)\rightarrow (t,r)$.
Similarly, we have
\begin{align} \label{ev6}
\begin{array}{l}
\|\tilde z_1-\tilde z_2\|_{C^{(1+\alpha)/2,1+\alpha}(\Delta^\infty_{T})} \\
\ \ \ \ \ \ \ \ \ \ \ \ \ \ \  \leq C_4(\|w_{1}-w_{2}\|_{C(\Delta_T)}+\|z_1-z_2\|_{L^\infty(\Delta^\infty_{T})}+\|h_1-h_2\|_{C^1([0,T])}),
\end{array}
\end{align}
where $C_4$ depends on $C_1, C_2$, the local Lipschitz coefficients of $f, g$ and the functions $A, B$ and $C$.
Taking the difference of the equations for $\overline h_1$ and $\overline
h_2$ results in
\begin{align}
\|\tilde h'_1-\tilde h'_2\|_{C^{\alpha/2}([0,T])}\leq
\mu\Big(\|\tilde w_{1, s}-\tilde w_{2, s}\|_{C^{\alpha/2,
0}(\Delta_T)}\Big). \label{e8}
\end{align}
 Combining \eqref{Bc}, (\ref{e6}), (\ref{ev6}) and (\ref{e8}), and assuming $T\leq 1$ we obtain
\begin{align*}
&\|\tilde w_1-\tilde w_2\|_{C^{(1+\alpha)/2,
1+\alpha}(\Delta_T)}+\|\tilde z_1-\tilde z_2\|_{C^{(1+\alpha)/2,
1+\alpha}(\Delta^\infty_T)}+\|\tilde h'_1-\tilde h'_2\|_{C^{\alpha/2}([0,T])}\\
&\quad \leq C_5
(\|w_{1}-w_{2}\|_{C(\Delta_T)}+\|z_1-z_2\|_{L^\infty(\Delta^\infty_T)}+\|h'_1-h'_2\|_{C[0,T]}),
\end{align*}
with $C_5$ depending on $C_3, C_4$ and $\mu$.  Hence for
\begin{align*}
T:= \min \left\{1,\,
\left(\frac 1{2C_5}\right)^{2/\alpha},\, (\mu C_1)^{-2/\alpha},\ C_1^{-2/(1+\alpha)}, \, \frac{h_0}{8(1+h^*)} \right\},
\end{align*}
we have
\begin{align*}
&\|\tilde w_1-\tilde w_2\|_{C(\Delta_T)}+\|\tilde z_1-\tilde z_2\|_{L^\infty(\Delta^\infty_T)}+\|\overline h'_1-\overline h'_2\|_{C([0,T])}\\
&\leq T^{(1+\alpha)/2}(\|\tilde w_1-\tilde w_2\|_{C^{(1+\alpha)/2, 1+\alpha}(\Delta_T)}+\|\tilde z_1-\tilde z_2\|_{C^{(1+\alpha)/2, 1+\alpha}(\Delta^\infty_T)})+
T^{\alpha/2}\|\overline h'_1-\overline h'_2\|_{C^{\alpha/2}([0,T])}\\
&\leq C_5 T^{\alpha/2}(\|w_{1}-w_{2}\|_{C(\Delta_T)}+\|z_1-z_2\|_{L^\infty(\Delta^\infty_T)}+\|h'_1-h'_2\|_{C([0,T])})\\
&\leq \frac12(\|w_{1}-w_{2}\|_{C(\Delta_T)}+\|\tilde z_1-\tilde z_2\|_{L^\infty(\Delta^\infty_T)}+\|h'_1-h'_2\|_{C([0,T])}).
\end{align*}
This shows that for this $T$, $\mathcal{F}$ is a contraction mapping in $\Gamma_{T}$. It follows from
 the contraction mapping theorem that $\mathcal{F}$ has a unique fixed point $(w,z,h)$ in $\Gamma_T$. In
other words, $(w(t, s), z(t, s); h(t))$ is the solution of the
problem (\ref{Bb}) and therefore $(u(t, r), v(t, r); h(t))$ is the
solution of the problem (\ref{1f1}). Moreover, by using the Schauder
estimates, we have additional regularity of the solution, $h(t)\in C^{1+\alpha/2}[0,T]$,  $u\in C^{1+\alpha/2, 2+\alpha}(G_T)$ and $v\in
C^{1+\alpha/2, 2+\alpha}((0, T]\times (0, +\infty))$. Thus $(u(t, s), v(t, s); h(t))$ is the classical solution of the problem (\ref{1f1}),
where $G_T:=\{(t,r)\in\mathbb{R}^2:t\in (0, T], r\in (0,h(t))\}$. \epf

\begin{rmk} {\rm By a bounded solution $(u,v,h)$ we mean that  there exists $M_T$ such that $|u|\leq M_T$ and
$|v|\leq M_T$ in $[0, T]\times [0, \infty)$. We cannot confirm the uniqueness of the solution to (\ref{1f1})
without the assumption of boundedness since
$v$ is defined in an unbounded domain. For our problem \eqref{f1},
the solution is always bounded, see Theorem \ref{globoun2} below. }
\end{rmk}

\begin{rmk}{\rm
It follows from the uniqueness of the solution to \eqref{1f1} and
a standard compactness argument that the unique solution $(u, v, h)$
depends continuously on the parameters appearing in \eqref{1f1}. This
fact will be used in the sections below.}
\end{rmk}

\begin{thm}\label{globoun}  Under the assumptions of Theorem \ref{localex}, if we assume further that
 there exists a constant $L>0$ such that $f(u, v)\leq L (u+v)$ and $g(u,v)\leq L (u+v)$ for $u, v\geq 0$, then the unique solution
 obtained in Theorem \ref{localex} can be extended uniquely to all $t>0$.
\end{thm}
\bpf  Let $[0, T_{max})$ be the maximal time interval in which the
solution exists. By Theorem \ref{localex},  $T_{max}>0$. It remains to show
that $T_{max}=\infty$.

Suppose for contradiction that $T_{max}<\infty$. Fix $M^*\in (T_{max}, \infty)$.
Let $(U(t), V(t))$ be the solution to the following ODE system:
\begin{align}
\left\{
\begin{array}{lll}
U_{t}=L(U+V),\ V_t=L (U+V),&t>0, \\
U(0)=\|u_0\|_{C([0, h_0])},\ V(0)=\|v_0\|_{L^\infty([0, \infty))}.
\end{array} \right.\label{ode1}
\end{align}
It is easy to see that
\begin{align*}
0<U+V<(\|u_0\|_{C[0, h_0]}+\|v_0\|_{L^\infty([0, +\infty))})e^{2Lt}\leq (\|u_0\|_{C[0, h_0]}+\|v_0\|_{L^\infty([0, +\infty))})e^{2LM^*}
\end{align*}
for $t\in [0, T_{\max})$. Recalling the assumption that $f(u, v)\leq L (u+v)$ and $g(u,v)\leq L (u+v)$ and comparing $(u,v)$ with $(U, V)$ yield that for $t\in [0, T_{\max})$ and
$r\in [0, h(t)]$,
\begin{align*}
0<u(t, r)+v(t,r)\leq U(t)+V(t)\leq (\|u_0\|_{C[0, h_0]}+\|v_0\|_{L^\infty([0, +\infty))})e^{2LM^*}:=C_6.
\end{align*}

 Next we claim that $0<h'(t)\leq C_7$ for all $t\in (0,T_{\max})$ and
some $C_7$ independent of $T_{\max}$. In fact, by the strong maximum principle and Hopf boundary lemma
$h'(t)$ is always positive as long as the solution exists.
To derive an upper bound of $h'(t)$, we define
\begin{align*}
\Omega=\Omega_M: =\{(t, r): 0<t<T_{\max},\, \, h(t)-M^{-1}<r<h(t)\}
\end{align*} and
construct an auxiliary function
\begin{align*}
\overline u(t, r):=C_6[2M(h(t)-r)-M^2(h(t)-r)^2].
\end{align*}
We will choose $M$ so that $\overline u(t,r)\geq u(t,r)$ holds over $\Omega$.

Direct calculations show that, for $(t,r)\in\Omega$,
\begin{align*}
&\overline u_t=2C_6Mh'(t)(1-M(h(t)-r))\geq 0,\\
&-\overline u_{r}=2C_6M[1-M(h(t)-r)]\geq 0,\\
&-\Delta \overline u=-\overline u_{rr}-\frac{N-1}r \overline u_r\geq 2C_6M^2,\\
&f(u, v)\leq L(u+v)\leq LC_6.
\end{align*}
It follows that
\begin{align*}
\overline u_t-d_1 \Delta \overline u\geq 2 d_1 C_6M^2\geq f(u,v) \mbox{  in }
\Omega
\end{align*} if $M^2\geq \frac {L}{2d_1}$. On the other hand,
\begin{align*}\overline u(t, h(t)-M^{-1})=C_6\geq u(t, h(t)-M^{-1}), \quad
\overline u(t, h(t))=0=u(t, h(t)).
\end{align*}

To use the maximum principle over $\Omega$, we only have to find some $M$ independent of
$T_{\max}$ such that $u_0(r)\leq \overline u(0,r)$ for $r\in [h_0-M^{-1}, h_0]$. We
calculate
\begin{align*}
\overline u_r(0,r)=-2C_6M[1-M(h_0-r)]\leq -C_6M \mbox{ for } r\in
[h_0-(2M)^{-1}, h_0].
\end{align*}
Therefore upon choosing
\begin{align*}
M:=\max\left\{ \sqrt{\frac{L}{2d_1}},\
\frac{4\|u_0\|_{C^1([0,h_0])}}{3C_6}\right\},
\end{align*}
we will have
\begin{align*}
\overline u_r(0,r)\leq u_0'(r) \mbox{ for } r\in
[h_0-(2M)^{-1}, h_0].
\end{align*}
Since $\overline u(0,h_0)=u_0(h_0)=0$, the above inequality implies
\begin{align*}
\overline u(0,r)\geq u_0(r) \mbox{ for } r\in
[h_0-(2M)^{-1}, h_0].
\end{align*}
Moreover, for $r\in [h_0-M^{-1}, h_0-(2M)^{-1}]$, we have
\begin{align*}
\overline u(0,r)\geq \frac 3 4 C_6,\; u_0(r)\leq
\|u_0\|_{C^1([0,h_0])}M^{-1}\leq \frac 34 C_6.
\end{align*}
Therefore $u_0(r)\leq \overline u(0,r)$ for $r\in [h_0-M^{-1}, h_0]$.

Applying the maximum principle to $\overline u-u$ over $\Omega$ gives that
$u(t,r)\leq \overline u(t,r)$ for $(t,r)\in\Omega$, which implies that
\[u_r(t, h(t))\geq \overline u_r(t, h(t))=-2MC_6,\; h'(t)=-\mu u_r(t,
h(t))\leq C_7:=2MC_6\mu\]
for $t\in [0, T_{\max})$. Moreover,
\begin{align*}
h_0\leq h(t)\leq h_0+C_7t\leq h_0+C_7M^*.
\end{align*}

We now fix $\delta_0\in (0, T_{max})$. By standard
parabolic regularity, we can find $C_8>0$ depending only on
$M^*$, $L, C_6$ and $C_7$ such that
 $\|u(t,\cdot)\|_{C^{1+\alpha}([0, h(t)])}\leq
C_8$ and $\|v(t,\cdot)\|_{C^{1+\alpha}([0, \infty))}\leq
C_8$ for $t\in [\delta_0, T_{\max})$.  It then follows from the
proof of Theorem \ref{localex} that there exists a $\tau>0$ depending only on
$M^*, L, C_6, C_7$ and $C_8$ such that the solution of problem \eqref{1f1}
with initial time $T_{max}-\tau/2$ can be extended uniquely to the
time $T_{max}-\tau /2+\tau$. This contradicts the maximality of $T_{max}$.
\epf

We have the following estimates.

\begin{thm}\label{globoun2} Problem \eqref{f1} admits a unique and uniformly bounded solution $(u,v,h)$. That is, the solution is defined for all $t>0$ and 
there exist constants $M_1$ and $M_2$ such
that
\begin{align*}
&0<u(t,r)\leq M_1\; \mbox{ for } t\in (0,+\infty),\, 0\leq r<h(t), \\
&0<v(t,x)\leq M_2\; \mbox{ for } t\in (0,+\infty),\, 0\leq r<+\infty.
\end{align*}
Moreover, there exist a constant $M_3$ such that
$$0<h'(t)\leq M_3 \; \mbox{ for }  t\in (0,+\infty).$$
Further more, \eqref{f1} does not have any unbounded solution.
\end{thm}
\bpf  By Theorem \ref{globoun}, \eqref{f1} has a unique bounded solution defined for all $t>0$.
It follows from the comparison principle that $u(t,
r)\leq \overline u(t)$ for $t\in(0, \infty)$ and $r\in [0, h(t)]$,
where
\begin{align*}
\overline u(t) :=\frac{a_1}{b_1} e^{\frac {a_1}{b_1} t}\Big(e^{\frac {a_1}{b_1} t}-1+
\frac {a_1}{b_1\|u_0\|_\infty}\Big)^{-1},
\end{align*}
which is the solution of the problem
\begin{align}
\label{ode}
\left\{
\begin{array}{l}
\displaystyle\frac {d \overline u}{d t}= \overline u (a_1-b_1\overline
u),\quad t>0, \\
  \overline u(0)=\|u_0\|_\infty.
  \end{array}\right.
\end{align}
Thus we have
\begin{align*}
u(t,r)\leq M_1:=\sup_{t\geq 0}\overline u(t).
\end{align*}

Since $v(t, r)$ satisfies
\begin{align*}
\left\{
\begin{array}{lll}
v_{t}-d_2\Delta v\leq v(a_2-c_2v),\; & t>0, \, 0\leq r<\infty, \\
v(0, r)=v_{0}(r)\geq 0,\; &0\leq r<\infty.
\end{array} \right.
\end{align*}
we have $v(t, r)\leq \max \{ \|v_0\|_{L^\infty (0, +\infty)},
\frac{a_2}{c_2}\}\triangleq M_2$.

Using the strong maximum principle to the equation of $u$ we
immediately obtain
\begin{align*}
u(t,r)>0,\;\; u_r(t, h(t))<0 \ \;\; \textrm{for} \ t>0, 0\leq r<h(t).
\end{align*}
 Hence $h'(t)>0$ for $t\in (0, \infty)$. Similarly we have $v(t,r)>0$ for $t>0, 0\leq r<\infty.$

It remains to show that $h'(t)\leq M_3$ for $t\in (0,+\infty)$ and
some $M_3$. The proof is similar as that of Theorem \ref{globoun} with $C_6$ replaced by $M_1$ and $M_3=C_7=2MM_1\mu$, we omit the details.

We next show that any solution of \eqref{f1} is bounded, namely, 
there exists $M>0$ such that $u,\; v\leq M$ in the range they are defined,  whenever $(u,v,h)$ is a solution to \eqref{f1} defined in some maximal interval $t\in (0, T)$.
Indeed, let $U(x)$ be the unique boundary blow-up solution of
\[
-d_1\Delta U=U(a_1-b_1U) \mbox{ in } B_1(0):=\{r<1\},\; U=\infty \mbox{ on } \partial B_1(0),
\]
and denote $\tilde u(t,x)=u(t, |x|)$;
then it is easily checked by using the comparison principle that $\tilde u(t,x_0+x)\leq \|u_0\|_\infty+U(x)$ for $x\in B_1(x_0)$ and $t>0$.
It follows that $u\leq \|u_0\|_\infty+U(0)$ in the range that $u$ is defined. Similarly we can show $v\leq \|v_0\|_\infty+V(0)$,
where $V(x)$ is the unique boundary blow-up solution of 
\[
-d_2 \Delta V=V(a_2-c_2V) \mbox{ in } B_1(0),\; V=\infty \mbox{ on } \partial B_1(0).
\]
(The existence and uniqueness  of $U$ and $V$ is well known; see, for example, \cite{DM}.)
 \epf

In what follows, we discuss the comparison principle for \eqref{f1}.
For a given pair of functions $\underline{\bm{u}}:=(\underline u, \underline v)$ and $\overline{\bm{u}}:=(\overline u, \overline v)$, we denote
\begin{align*}
[\underline{\bm{u}}, \overline{\bm{u}}] = \{ \bm{u} :=(u, v)\in [C([0, T]\times [0, \infty))]^2: \ (\underline u, \underline v)
\leq (u,v)\leq (\ol u, \ol v)\},
\end{align*}
where by $(u_1, v_1)\leq (u_2, v_2)$, we mean $u_1\leq u_2$ and $v_1\leq v_2$.
A function pair $(f, g)=(f(u,v), g(u,v))$ is said to be quasimonotone nonincreasing if for fixed $u$, $f$ is nonincreasing in $v$,
and for fixed $v$, $g$ is nonincreasing in $u$; this is satisfied by $f=u(a_1-b_1u-c_1v)$ and $g=v(a_2-b_2u-c_2v)$ in \eqref{f1} for $u, v\geq 0$.

\begin{lem}[The Comparison Principle]\label{comparison} Let $(f,g)$ be quasimonotone nonincreasing and Lipschitz continuous in
$[\underline{\bm{u}}, \overline{\bm{u}}]$, with $f(0,v)=g(u,0)\equiv 0$.
  Assume that $T\in (0,\infty)$, $\underline h, \overline
h\in C^1([0,T])$, $\underline u\in C(\overline{D^*_T})\cap
C^{1,2}(D^*_T)$ with $D^*_T:=\{(t,r)\in\mathbb{R}^2:t\in(0,T], r\in(0,\underline{h}(t))\}$, $\overline u\in C(\overline{D^{**}_T})\cap
C^{1,2}(D^{**}_T)$ with $D^{**}_T:=\{(t,r)\in\mathbb{R}^2:t\in(0,T], r\in[0,\overline{h}(t))\}$, $\underline v, \overline v\in
(L^\infty\cap C)([0, T]\times [0, \infty))\cap
C^{1,2}((0, T]\times [0, \infty))$ and
\begin{eqnarray*}
\left\{
\begin{array}{lll}
\overline u_{t}-d_1 \Delta \overline u\geq f(\overline u, \underline v),\; & 0<t\leq T, \ 0\leq r<\overline h(t), \\
\underline u_{t}-d_1 \Delta \underline u\leq f(\underline u, \overline v),\; & 0<t\leq T, \ 0\leq r<\underline h(t), \\
\overline v_{t}-d_2 \Delta \overline v\geq g(\underline u, \overline v),\; & 0<t\leq T, \ 0\leq r<\infty,  \\
\underline v_{t}-d_2 \Delta \underline v\leq g(\overline u, \underline v),\; & 0<t\leq T, \ 0\leq r<\infty,  \\
\overline u_r(t,0)=\underline v_r(t,0)=0,\; \overline  u(t, r)=0,\quad & 0<t\leq T,\ \overline h(t)\leq r<\infty,\\
\underline u_r(t,0)=\overline v_r(t,0)=0,\; \underline  u(t, r)=0,\quad & 0<t\leq T,\ \underline h(t)\leq r<\infty,\\
\underline h'(t)\leq -\mu \underline u_r(t,h(t)),\, \overline h'(t)\geq -\mu \overline u_r(t,h(t)),\quad &0<t\leq T,\\
\underline h(0)\leq h_0\leq \overline h(0),   &\\
\underline u(0, r)\leq u_{0}(r)\leq \overline u(0,r),\;  &0\leq r\leq h_0,\\
\underline v(0, r)\leq v_0(r)\leq \overline v(0, r), &0\leq r<\infty.
\end{array} \right.
\end{eqnarray*}
Let $(u,v,h)$ be the unique bounded solution of \eqref{1f1}. Then
\begin{align*}
&h(t)\le\overline h(t)\ {\it in}\ (0, T],\ u(t, r)\leq\overline u(t, r),\ v(t, r)\geq \underline v(t, r)\ {\it for}\ (t, r)\in (0, T]\times[0, \infty), \\
&h(t)\geq\underline h(t)\ {\it in}\ (0, T],\ u(t, r)\geq\underline u(t, r),\ v(t,r)\leq \overline v(t,r)\ {\it for}\ (t, r)\in (0, T]\times[0, \infty).
\end{align*}

\end{lem}
\bpf
We only  prove  $u\leq \overline u$, $v\geq \underline v$ and $h\leq \overline h$; the result involving $(\underline u, \overline v, \underline h)$
can be proved in a similar way.
Let $\tilde M$ be an upper bound of $v$ and $\underline v$ in $[0, T]\times [0, +\infty)$,
$w=\tilde M-v$ and $\overline w=\tilde M-\underline v$, then $(\overline u, \overline w, \overline h)$ satisfies
\begin{align}
\left\{
\begin{array}{lll}
\overline u_{t}-d_1 \Delta \overline u\geq f(\overline u, \tilde M-\overline w),\; & 0<t\leq T, \ 0\leq r<\overline h(t), \\
\overline w_{t}-d_2 \Delta \overline w\geq -g(\overline u, \tilde M-\overline w),\; & 0<t\leq T, \ 0\leq r<\infty,  \\
\overline u_r(t,0)=\overline w_r(t,0)=0,\; \overline  u(t, r)=0,\quad &0<t\leq T,\ \overline h(t)\leq r<\infty,\\
-\mu \overline u_r(t,h(t))\leq \overline h'(t),\quad &0<t\leq T,\\
h_0\leq h(0), \ u_{0}(r)\leq \overline u(0,r),\;  &0\leq r\leq \overline h_0,\\
\tilde M-v_0(r)\leq \overline w(0, r), &0\leq r<\infty.
\end{array} \right.\label{cp}
\end{align}

First assume that $h_0<\overline h(0)$. We claim that $h(t)<\overline h(t)$ for all $t\in (0, T]$.
If our claim does not hold, then we
can find a first $t^*\leq T$ such that $h(t)<\overline
h(t)$ for $t\in (0, t^*)$ and $h(t^*)=\overline h(t^*)$. It
follows that
\begin{align}
\label{h-h}
 h'(t^*)\geq \overline h'(t^*).
\end{align}
 We now show that $(u, w)\leq (\overline u, \overline w)$ in $[0, t^*]\times [0, \infty)$.
 Letting $U=(\overline u-u)e^{-Kt}$ and $W=(\overline w-w)e^{-Kt}$, we obtain
 \begin{align}
\left\{
\begin{array}{ll}
U_{t}-d_1 \Delta U \geq (-K+b_{11}) U+b_{12}W,\; & 0<t\leq t^*, \ 0\leq r<h(t), \\
W_{t}-d_2 \Delta W\geq b_{21}U+(-K+b_{22})W,\; & 0<t\leq t^*, \ 0\leq r<\infty,  \\
U_r(t,0)=W_r(t,0)=0,\; U(t, r)=0,\quad &0<t\leq t^*,\ \overline{h}(t)\leq r<\infty,\\
U(0,r)\geq 0, W(0,r)\geq 0, &0\leq r<\infty,
\end{array} \right.\label{cp1}
\end{align}
 where $K$ is sufficiently large such that $K\geq 1+|b_{11}|+b_{12}+|b_{22}|+b_{21}$ in $[0, t^*]\times [0, \infty)$,
  $b_{11}(t, r),  b_{22}(t, r)$ are bounded, $b_{21}(t,r), b_{12}(t, r)$ are bounded and nonnegative
  since  $(f,g)$ in (\ref{cp}) is quasimonotone nonincreasing  and Lipschitz continuous.

 Since the first inequality of (\ref{cp1}) holds only in part of $[0,\infty)$,
 we cannot use the maximum principle directly.  We first prove that for any $l>h(t^*)$,
 \begin{align*}
 U(t, r)\geq - \frac {\tilde M(r^2+2dN t)}{l^2}\ \textrm{and}\
   W(t, r)\geq -\frac {\tilde M(r^2+2dN t)}{l^2}
   \end{align*}
 in $[0, t^*]\times [0,l]$, where $d=\max(d_1, d_2)$. (Recall that $\tilde M$ is an upper bound of $v$ and $\underline v$.)

We observe that due to the inequalities satisfied by $\overline u$, we can apply the maximum principle  to $\overline u$
over the region $\{(t,r): 0\leq r\leq \overline h(t), 0\leq t\leq T\}$ to conclude that $\overline u\geq 0$.
 Set 
 \[
 \mbox{ $\ol U(t, r)= U+\frac {\tilde M(r^2+2dN t)}{l^2}$ and
 $\ol W(t, x)=W+\frac {\tilde M(r^2+2d Nt)}{l^2}$;}
 \]
 then  due to our choice of $K$, $(\ol U, \ol W)$ satisfies
 \begin{align*}
\left\{
\begin{array}{ll}
\ol U_{t}-d_1 \Delta \ol U \geq (-K+b_{11}) \ol U+b_{12}\ol W,\; & 0<t\leq t^*, \ 0\leq r<h(t), \\
\ol W_{t}-d_2 \Delta \ol W\geq b_{21}\ol U+(-K+b_{22})\ol W,\; & 0<t\leq t^*, \ 0\leq r<l,  \\
\ol U_r(t,0)=\ol W_r(t,0)=0,\; \ol U(t, r)\geq \frac {\tilde M(r^2+2d Nt)}{l^2}>0,\quad &0<t\leq t^*,\ h(t)\leq r \leq l,\\
\ol W(t, l)=W(t,l)+\frac {\tilde M(l^2+2d Nt)}{l^2}>0, &0<t\leq t^*,\\
\ol U(0,r)\geq 0, \ol W(0,r)\geq 0, &0\leq r\leq l.
\end{array} \right.
\end{align*}
We now prove that $\min \{\min_{[0, t^*]\times [0, l]} \ol U, \min_{[0, t^*]\times [0, l]} \ol W\}:=\tau\geq 0$.
In fact, if $\tau<0$, then there exists $(t_0,r_0)\in\mathbb{R}^2$ with $0<t_0\le t^*$ and $0\leq r_0<h(t_0)$ such that $\ol U(t_0, r_0)=\tau<0$, or there exists
$(t_1,r_1)\in\mathbb{R}^2$ with $0<t_1\le t^*$ and $0\leq r_1<l$ such that $\ol W(t_1, r_1)=\tau<0$. For the former case, $(\ol U_{t}-d_1 \Delta \ol U)(t_0, r_0)\leq 0$, but
\begin{align*}
[(-K+b_{11}) \ol U+b_{12}\ol W](t_0, r_0)\geq (-K+|b_{11}|)\tau +b_{12}\tau \geq -\tau >0.
\end{align*}
For the latter case, $(\ol W_{t}-d_2 \Delta \ol W)(t_1, r_1)\leq 0$, but
\begin{align*}
[b_{21}\ol U+(-K+b_{22})\ol W](t_1, r_1)\geq (-K+|b_{22}|)\tau +b_{21}\tau \geq -\tau >0.
\end{align*}
Both are impossible. Therefore $\tau\geq 0$, that is $\ol U\geq 0$ and $\ol W\geq 0$ in $[0, t^*]\times [0, l]$, which implies that
\begin{align*}
U(t, r)\geq - \frac {\tilde M(r^2+2dN t)}{l^2}, \quad
   W(t, r)\geq -\frac {\tilde M(r^2+2dN t)}{l^2}
\end{align*}
  for $0\leq t\leq t^*$, $ 0\leq r\leq l$.
  Taking $l\to \infty$ yields that $U(t, r)\geq 0$ and $W(t, r)\geq 0$ in $[0, t^*]\times [0, \infty)$,
  therefore $u\leq \ol u$ and $w\leq \ol w$ in $[0, t^*]\times [0, \infty)$.

 We now compare
$u$ and $\overline u$ over the bounded region
\begin{align*}
\Omega_{t^*}:=\{(t, r)\in\R^2: 0< t\leq t^*, 0\leq r< h(t)\}.
\end{align*}
Since $Z(t,r):=\overline u(t,r)-u(t,r)$ satisfies
\begin{align*}
Z_{t}-d_1 \Delta Z \geq b_{11} Z+b_{12}(\ol w-w)\geq b_{11} Z,\; 0<t\leq t^*, \ 0\leq r<h(t),
\end{align*}
 the strong maximum principle and the Hopf boundary lemma yield
$Z(t,r)>0$ in $\Omega_{t^*}$, and $Z_r(t^*,h(t^*))<0$. We then deduce that $h'(t^*)<\overline
h'(t^*)$. But this contradicts \eqref{h-h}. This proves our claim
that $h(t)<\overline h(t)$ for all $t\in (0, T]$. We may
now apply the above procedure over $[0, T]\times [0, \infty)$ to
conclude that $u\leq \overline u$ and $w\leq \ol w$ (i.e. $v\geq \underline v$) in $[0, T]\times [0, \infty)$. Moreover,
$u<\overline u$ for $t\in(0, T]$ and $r\in [0, h(t))$.

If $h_0=\overline h(0)$, we use approximation. For small $\epsilon>0$, let $(u_\epsilon, v_\epsilon, h_\epsilon)$ denote
the unique solution of \eqref{f1} with $h_0$ replaced by
$h_0(1-\epsilon)$. Since the unique solution of \eqref{f1} depends continuously on the
parameters in \eqref{f1}, as $\epsilon\to 0$, $(u_\epsilon, v_\epsilon,
h_\epsilon)$ converges to $(u,v,h)$, the unique solution of
\eqref{f1}. The desired result then follows by letting $\epsilon\to
0$ in the inequalities  $u_\epsilon\leq\overline u$, $v_\epsilon\geq \underline v$ and
$h_\epsilon<\overline h$.
\epf

\begin{rmk}\label{condition-r=0}
{\rm The conclusions in Lemma \ref{comparison} remain valid if the condition
\[
z_r(t,0)=0 \mbox{ for } z\in\{\underline u, \overline u,\underline v,\overline v\}
\]
is replaced by
\[
z_r(t,0)>0 \mbox{ for } z\in \{\underline u,\underline v\},\; z_r(t,0)<0 \mbox{ for } z\in\{\overline u,\overline v\}.
\]
To see this, we only need to observe that  all the arguments in the proof carry over to the new case,
except that we now have to avoid $r_0=0$ or $r_1=0$, since the functions $\overline U(t,|x|)$ and $\overline W(t, |x|)$ are not
 $C^2$ in $x$ at $x=0$. However, the above conditions guarantee that
$\overline U_r(t,0)<0$ and $\overline W_r(t,0)<0$. Therefore $(t,0)$ cannot be a minimum point of these functions. This implies that
$r_0>0$ and $r_1>0$. }
\end{rmk}
We next fix $u_0, v_0, d_i, a_i, b_i, c_i$ and examine the dependence of the solution on $\mu$,
and we write $(u^{\mu}, v^{\mu}, h^{\mu})$ to emphasize this dependence. As a consequence of Lemma 2.6, we have the following result.

\begin{cor}\label{monotone} For fixed $u_0, v_0, d_i, a_i, b_i$ and $c_i$.
If $\mu_1\leq \mu_2$. Then $u^{\mu_1}(t, r)\leq u^{\mu_2}(t,r)$, $v^{\mu_1}(t, r)\geq v^{\mu_2}(t,r)$ for
 $t\in(0, \infty)$, $r\in [0, h^{\mu_1}(t))$ and $h^{\mu_1}(t)\leq h^{\mu_2}(t)$ in $(0, \infty)$.
\end{cor}

\section{Invasion of an inferior competitor}

In this section, we examine the case that $u$ is an inferior competitor, namely
\begin{equation}
\label{u-inf}
\frac{a_1}{a_2}<\min\left\{\frac{b_1}{b_2},\frac{c_1}{c_2}\right\}.
\end{equation}

The following theorem shows that the inferior invader cannot establish itself and the native species always survives the invasion.

\begin{thm}\label{vanofu} If \eqref{u-inf} holds  and $v_0\not\equiv 0$, then
\[
\lim_{t\to +\infty} \ \big(u(t, r), v(t, r)\big)=\left(0,\,\frac
{a_2}{c_2}\right) \mbox{ uniformly in any compact subset of $[0, \infty)$.}
\]
\end{thm}
\bpf
First we recall that the comparison principle gives $u(t,
r)\leq u^*(t)$ for $t>0$ and $r\in [0, h(t)]$, where
\begin{align*}u^*(t) =\frac{a_1}{b_1} e^{\frac {a_1}{b_1} t}\Big(e^{\frac {a_1}{b_1} t}-1+\frac {a_1}{b_1\|u_0\|_\infty}\Big)^{-1}
\end{align*}
is the solution of the problem
\begin{align}
\label{ode1}
\left\{
\begin{array}{ll}
 (u^*)'= u^* (a_1-b_1 u^*),&\ t>0, \\
u^*(0)=\|u_0\|_\infty. &
\end{array}\right.
\end{align}
Since $\lim_{t\to\infty}
u^*(t)= \frac {a_1}{b_1}$, we deduce
\begin{align*}
\limsup_{t\to +\infty}u(t, r)\leq \frac {a_1}{b_1}\  \ \textrm{uniformly for}\ \ r\in [0,\infty).
\end{align*}
Similarly, we have
\begin{align}
\label{od1}
\limsup_{t\to +\infty}v(t, r)\leq \frac {a_2}{c_2}\  \ \textrm{uniformly for}\ \ r\in [0,\infty).
\end{align}
Therefore for $\varepsilon_1=(\frac {a_2}{b_2}-\frac {a_1}{b_1})/2$, there exists $t_1>0$ such that
$u(t, r)\leq \frac {a_1}{b_1}+\varepsilon_1$ for
$t\geq t_1, r\in [0,\infty)$. Then $v$ satisfies
  \begin{align}
\left\{
\begin{array}{lll}
v_{t}-d_2 \Delta v\geq v(b_2\varepsilon_1 -c_2 v),\; & t>t_1, \ 0\leq r<\infty,  \\
v_r(t,0)=0,\quad & t>t_1\\
v(t_1, r)>0, &0\leq r<\infty.
\end{array} \right.
\label{fv1}
\end{align}

Let $v_*$ be the unique solution to
\begin{align*}
\left\{
\begin{array}{lll}
(v_*)_{t}-d_2 \Delta v_*=v_*(b_2\varepsilon_1 -c_2 v_*),\; & t>t_1, \ 0\leq r<\infty,  \\
(v_*)_r(t,0)=0,\quad & t>t_1\\
v_*(t_1, r)=v(t_1, r), &0\leq r<\infty.
\end{array} \right.
\end{align*}
It is well known (\cite{DM}) that $\lim_{t\to \infty}v_*(t, r)=b_2\varepsilon_1/c_2$ uniformly in any bounded subset of $[0, \infty)$.
 Therefore for any $L>0$, there exists $t_L>t_1$ such that
\begin{align}
\label{ine}
v(t, r)\geq v_*(t, r)\geq \frac{b_2\varepsilon_1}{2c_2}\  \textrm{ for}\ t\geq t_L,\ 0\leq r\leq L.
\end{align}

Now $(u,v)$ satisfies
\begin{align}
\left\{
\begin{array}{lll}
u_{t}-d_1 \Delta u=u(a_1-b_1 u-c_1 v),\; & t>t_L, \ 0<r<h(t),  \\
v_{t}-d_2 \Delta v=v(a_2-b_2 u-c_2 v),\; & t>t_L, \ 0<r<\infty,  \\
u_r(t,0)=v_r(t,0)=0,\quad & t>t_L,\\
u(t, r)\leq \frac {a_1}{b_1}+\varepsilon_1,\ v(t, r)\geq \frac{b_2\varepsilon_1}{2c_2},  & t\geq t_L,\ 0\leq r\leq L.
\end{array} \right.
\label{fs1}
\end{align}
Since $u\equiv 0$ for $t>t_L, \ r\geq h(t)$, no matter whether or not $h(t)\leq L$, we always have $u\leq \overline u$
and $v\geq \underline v$ in $[t_L, \infty)\times [0, L]$,
where $(\overline u, \underline v)$ satisfies
\begin{align}
\left\{
\begin{array}{lll}
\overline u_{t}-d_1 \Delta \overline u=\overline u(a_1-b_1\overline u-c_1 \underline v),\; & t>t_L, \ 0<r<L,  \\
\underline v_{t}-d_2 \Delta \underline v=\underline v(a_2-b_2 \overline u-c_2 \underline v),\; & t>t_L, \ 0<r<L,  \\
\overline u_r(t,0)=\underline v_r(t,0)=0,\quad & t>t_L,\\
\overline u(t, r)=\frac {a_1}{b_1}+\varepsilon_1,\ \underline v(t, r)=\frac{b_2\varepsilon_1}{2c_2},
 & t\geq t_L,\ r=L, \textrm{or}\, t=t_L,\, 0\leq r\leq L.
\end{array} \right.
\label{fs11}
\end{align}
The system (\ref{fs11}) is quasimonotone nonincreasing, which generates a monotone dynamical system with respect to the order
\[
(u_1,v_1)\leq_P (u_2, v_2) \mbox{ if and only if } u_1\leq u_2\mbox{ and } v_1\geq v_2,
\]
with the initial value $\left(\frac {a_1}{b_1}+\varepsilon_1,\ \frac{b_2\varepsilon_1}{2c_2}\right)$ an upper solution.
It follows from the theory of monotone dynamical systems (see, e.g. \cite{hS} Corollary 3.6) that
$\lim_{t\to +\infty} \ \overline u(t, r)=\overline u_L(r)$ and $\lim_{t\to +\infty} \ \underline v(t, r)=\underline v_L(r)$ uniformly in
$[0, L]$, where $(\overline u_L, \underline v_L)$  satisfies
\begin{align} \label{fs12}\left\{
\begin{array}{ll}
-d_1 \Delta \overline u_L=\overline u_L(a_1-b_1 \overline u_L-c_1\underline v_L),\; &\ 0\leq r<L,  \\
-d_2 \Delta \underline v_L=\underline v_L(a_2-b_2 \overline u_L-c_2\underline v_L),\; &\ 0\leq r<L, \\
\frac{\partial \overline u_L}{\partial r}(0)=\frac{\partial \underline v_L}{\partial r}(0)=0,\; & \\
\overline u_L(L)=\frac {a_1}{b_1}+\varepsilon_1,\ \underline v_L(L)=\frac{b_2\varepsilon_1}{2c_2}, &
\end{array} \right.
\end{align}
and is the maximal solution below $(\frac{a_1}{b_1}+\varepsilon_1, \frac{b_2\varepsilon_1}{2c_2})$ of the above problem under the order $\leq_P$.

Next we observe that if $0<L_1<L_2$, then $\overline u_{L_1}(r)\geq \overline u_{L_2}(r)$ and $\underline v_{L_1}(r)\leq \underline v_{L_2}(r)$ in $[0, L_1]$, which can be derived by comparing the boundary conditions and initial conditions in (\ref{fs11}) for $L=L_1$ and $L=L_2$.

Let $L\to \infty$, by classical elliptic regularity theory and a diagonal procedure, it follows that  $(\overline u_{L}(r), \underline v_{L}(r))$ converges
uniformly on any compact subset of $[0, \infty)$ to $(\overline u_\infty, \underline v_\infty)$, which satisfies
\begin{align*} \left\{
\begin{array}{lll}
-d_1 \Delta \overline u_\infty=\overline u_\infty(a_1-b_1 \overline u_\infty-c_1\underline v_\infty),\; &\ 0\leq r<\infty,  \\
-d_2 \Delta \underline v_\infty=\underline v_\infty(a_2-b_2 \overline u_\infty-c_2\underline v_\infty),\; &\ 0\leq r<\infty, \\
\frac{\partial \overline u_\infty}{\partial r}(0)=\frac{\partial \underline v_\infty}{\partial r}(0)=0,\; & \\
\overline u_\infty(r)\leq\frac {a_1}{b_1}+\varepsilon_1,\ \underline v_\infty(r)\geq\frac{b_2\varepsilon_1}{2c_2}, &0\leq r<\infty.
\end{array} \right.
\end{align*}

Next we show that $\overline u_\infty(r)\equiv 0$ and $\underline v_\infty(r)\equiv \frac{a_2}{c_2}$. To this end, let us consider the following
ODE system:
\begin{align} \left\{
\begin{array}{lll}
z_t=z(a_1-b_1 z-c_1 w),\; &t>0,  \\
w_t=w(a_2-b_2 z-c_2 w),\; &t>0, \\
z(0)=\frac {a_1}{b_1}+\varepsilon_1,\; w(0)=\frac{b_2\varepsilon_1}{2c_2}.
\end{array} \right.
\label{gs1}
\end{align}
Since $a_1/a_2<\min \{c_1/c_2, b_1/b_2\}$, it is well-known (e.g.\cite{MT}) that $(z, w)\to (0, \frac {a_2}{c_2})$ as $t\to \infty$.
Therefore the solution $(Z(t, r), W(t, r))$ of the problem
\begin{align} \left\{
\begin{array}{lll}
Z_t-d_1\Delta Z=Z(a_1-b_1 Z-c_1 W),\; &t>0, \ r\geq 0, \\
W_t-d_2\Delta W=W(a_2-b_2 Z-c_2 W),\; &t>0, \ r\geq 0,\\
Z_r(t, 0)=W_r(t,0)=0,&t>0,\\
Z(0, r)=\frac {a_1}{b_1}+\varepsilon_1,\; W(0, r)=\frac{b_2\varepsilon_1}{2c_2}, & r\geq 0
\end{array} \right.
\label{gs2}
\end{align}
satisfies  $(Z, W)\to (0, \frac {a_2}{c_2})$ as $t\to \infty$ uniformly in $[0, \infty)$.
It follows from the comparison principle that $\overline u_\infty(r)\leq Z(t, r)$ and $\underline v_\infty(r)\geq  W(t, r)$ for $t>0$,
which immediately gives
that $\overline u_\infty=0$ and $\underline v_\infty=\frac {a_2}{c_2}$.

We thus have $\limsup_{t\to +\infty}u(t, r)\leq 0$ and $\liminf_{t\to +\infty}v(t, r)\geq \frac {a_2}{c_2}$ and
uniform in $[0,L]$, which together with (\ref{od1}) implies that
 $\lim_{t\to +\infty}u(t, r)=0$
and  $\lim_{t\to +\infty}v(t, r)=\frac
{a_2}{c_2}$ uniformly in any bounded subset of $[0, \infty)$.
\epf
\smallskip

\begin{rmk}
\label{rmk-u-inf}{\rm
The above proof can be used to show that the unique solution $(u(t,x), v(t,x))$
of the Cauchy problem of \eqref{entire}, with $u(0,x), v(0,x)$ nonnegative, bounded and $v(0,x)\not\equiv 0$ $($but not necessarily radially symmetric$)$,
also converges to $R_2=(0, a_2/c_2)$ locally uniformly in $\R^N$ as $t\to\infty$ when \eqref{u-inf} holds.}
\end{rmk}

We note that Theorem \ref{vanofu} gives no information on the dynamical behavior of the spreading front, and the exact behavior of $(u,v)$ over the entire spatial range $0\leq r<\infty$ is also unclear.
Our next result provides such information, provided that the native species is already rather established
at $t=0$, in the sense that $\inf_{r\geq 0}v_0(r)>0$.

 Let us observe that  Theorem \ref{globoun2} implies $h(t)$ is monotonic increasing and therefore
 there exists  $h_\infty\in (0, +\infty]$ such that $\lim_{t\to +\infty} \ h(t)=h_\infty$.

\begin{thm}\label{vanofu2} Suppose that \eqref{u-inf} holds and $v_0(r)\geq \delta >0$ for $0\leq r<\infty$. Then $h_\infty<\infty$ and $\lim_{t\to +\infty}u(t, r)=0$, $\lim_{t\to +\infty}v(t, r)=\frac
{a_2}{c_2}$ uniformly in $[0, \infty)$.
\end{thm}
\bpf
First it follows from the comparison principle that
\begin{align*}
\limsup_{t\to +\infty}u(t, r)\leq \frac {a_1}{b_1}\ \ \textrm{uniformly for}\ \ r\in [0,\infty). \\
\limsup_{t\to +\infty}v(t, r)\leq \frac {a_2}{c_2}\ \ \textrm{uniformly for}\ \ r\in [0,\infty).
\end{align*}
Therefore for $\varepsilon_1=(\frac {a_2}{b_2}-\frac {a_1}{b_1})/2$, there exists $t_1>0$ such that $u(t, r)\leq \frac {a_1}{b_1}+\varepsilon_1$ for
$t\geq t_1, r\in [0,\infty)$.

Note that $0\leq u(t, r)\leq M_1:=\max\{\frac {a_1}{b_1},\, \|u_0\|_{C([0, h_0])}\}$ and
$0\leq v(t, r)\leq M_2:=\max\{\frac {a_2}{c_2},\, \|v_0\|_{L^\infty([0, +\infty))}\}$ for $t>0, 0\leq r<\infty.$
Therefore $v$ satisfies
  \begin{align*}
\left\{
\begin{array}{ll}
v_{t}-d_2 \Delta v\geq v(a_2-b_2M_1-c_2 M_2),\; & t>0, \ 0\leq r<\infty,  \\
v_r(t,0)=0,\quad & t>0\\
v(0, r)\geq \delta, &0\leq r<\infty,\\
\end{array} \right.
\end{align*}
and hence $v(t, r)\geq \delta e^{(-b_2M_1-c_2 M_2)t}$ for $t>0$ and $0\leq r<\infty.$

Now let us consider the following problem:
\begin{align} \left\{
\begin{array}{lll}
z_t=z(a_1-b_1 z-c_1 w),\; &t>t_1,  \\
w_t=w(a_2-b_2 z-c_2 w),\; &t>t_1, \\
z(t_1)=\frac {a_1}{b_1}+\varepsilon_1,\; w(t_1)=\delta e^{(-b_2M_1-c_2 M_2)t_1}.
\end{array} \right.
\label{gf1}
\end{align}
It follows from the comparison principle that $u(t, r)\leq z(t)$ and $v(t, r)\geq w(t)$ for $t\geq t_1, 0\leq r <\infty$.
Under the assumption  $a_1/a_2<\min \{c_1/c_2, b_1/b_2\}$, it is well-known  that $(z, w)\to (0, \frac {a_2}{c_2})$ as $t\to \infty$.
It follows that
$\lim_{t\to +\infty}u(t, r)=0$ uniformly for $r\in[0,\infty)$ (we note that $u(t,r)=0$ for $r\ge h(t)$). Next we prove $\lim_{t\to\infty}v(t,r)=\frac{a_2}{c_2}$
uniformly for $r\in[0,\infty)$.
Since $\lim_{t\to\infty}u(t,r)=0$ uniformly for $r\in[0,\infty)$, for any $\varepsilon>0$ there exists $T>0$ such that $0\le u(t,r)\le \varepsilon$ for $t\ge T$ and
$r\in[0,\infty)$. Recall that  $v(t,T)\ge\delta e^{(-b_2M_1-c_2M_2)T}$. We thus have
\begin{align*}
\left\{
\begin{array}{ll}
v_t-d_2\Delta v\ge v(a_2-b_2\varepsilon-c_2v),&\ t\ge T,\ 0\le r< \infty, \\
v_r(t,0)=0, &\ t>T, \\
v(T,r)\ge\delta e^{(-b_2M_1-c_2M_2)T}, & 0\le r<\infty
\end{array}\right.
\end{align*}
 Consider the following problem
 \begin{align*}
\tilde{v}_t=\tilde{v}(a_2-b_2\varepsilon-c_2\tilde{v}) \mbox{ for } t\ge T, \;
\tilde{v}(T)=\delta e^{(-b_2M_1-c_2M_2)T}.
\end{align*}
It follows from the comparison principle that $\tilde{v}(t)\le v(t,r)$ for $t\ge T$ and $r\in[0,\infty)$. Since $\tilde{v}(t)\to\frac{a_2-\varepsilon b_2}{c_2}$ as $t\to\infty$, we have
$\frac{a_2-\varepsilon b_2}{c_2}\le\liminf_{t\to\infty}v(t,r)$ uniformly for $r\in[0,\infty)$. Since $\varepsilon>0$ is arbitrary we have $\lim_{t\to\infty}v(t,r)=\frac{a_2}{c_2}$ uniformly in $[0,\infty)$.

Next we prove that $h_\infty<\infty$.
Since $\lim_{t\to +\infty}u(t, r)=0$
and  $\lim_{t\to +\infty}v(t, r)=\frac
{a_2}{c_2}$ uniformly in $[0, \infty)$, we have $\lim_{t\to +\infty}(a_1-b_1u(t, r)-c_1v(t, r))=a_1-a_2c_1/c_2<0$
uniformly in $[0, \infty)$, and there exists $T^*$ such that
$(a_1-b_1u(t, r)-c_1v(t, r))\leq 0$ for $t\geq T^*, 0\leq r<\infty$.

Direct calculation gives
\begin{align*}
\frac{\textrm{d}}{\textrm{d} t}\int_0^{h(t)}r^{N-1}u(t, r)\textrm{d}r
&\ =\int_0^{h(t)}r^{N-1}u_t(t, r)\textrm{d}r+h^{N-1}(t)h'(t)u(t, h(t))\\
&\ =\int_0^{h(t)}d_1r^{N-1}\Delta u\textrm{d}r+\int_0^{h(t)}u(a_1-b_1 u-c_1v)r^{N-1}\textrm{d}r\\
&\ =-\frac{d_1}{\mu}h^{N-1}(t)h'(t)+\int_0^{h(t)}u(a_1-b_1 u-c_1v)r^{N-1}\textrm{d}r.
\end{align*}
Integrating from $T^*$ to $t$ yields
\begin{align}\label{k1}
\begin{array}{ll}
\displaystyle 0\leq \int_0^{h(t)}r^{N-1}u(t, r)\textrm{d}r =&\displaystyle \int ^{h(T^*)}_0 r^{N-1}u(T^*, r)dr +\frac {d_1}{N\mu}h^N(T^*)-\frac {d_1}{N\mu}h^N(t) \\
\displaystyle &\displaystyle+\int_{T^*}^t\int_0^{h(s)}u(a_1-b_1 u-c_1v)r^{N-1}drds, \quad t\geq T^*.
\end{array}
\end{align}
Since $(a_1-b_1u(t, r)-c_1v(t, r))\leq 0$ for $t\geq T^*, 0\leq r<\infty$, we  deduce
\begin{align*}
h^N(t)\leq \frac {N\mu}
{d_1}\int ^{h(T^*)}_0 r^{N-1}u(T^*,r)dr +h^N(T^*),\ t\geq T^*,
\end{align*}
which implies that $h_\infty<\infty$.
\epf

 Theorem \ref{vanofu2} suggests that an inferior competitor can never penetrate deep into the habitat of a well established native species,
 and it dies out before its invading front reaches a certain finite limiting position.

\section{Invasion of a superior competitor}

This section is devoted to the case that $u$ is a superior competitor, that is
\begin{equation}
\label{u-sup}
\frac{a_1}{a_2}>\max \left\{\frac{b_1}{b_2},\, \frac{c_1}{c_2}\right\}.
\end{equation}

Let $\lambda_1(R)$ be the principal eigenvalue of the operator $-\Delta$ in $B_R$ subject to homogeneous Dirichlet boundary conditions.
It is well-known that $\lambda_1(R)$ is a strictly decreasing continuous function and
\begin{align*}
\lim_{R\to 0^+}\lambda_1(R)=+\infty \ \textrm{and}\ \lim_{R\to +\infty}\lambda_1(R)=0.
\end{align*}
Therefore, there exists a unique $R^*$ such that
\begin{align}\label{lmd}
\lambda_1(R^*)=1.
\end{align}
It is easy to check that $R^*=\pi/2$ if $N=1$.

For the sake of convenience and completeness, we first recall the following
spreading-vanishing dichotomy for
the  radially symmetric diffusive logistic problem
\begin{align}
\left\{
\begin{array}{ll}
u_{t}-d \Delta u=u(a-b u),\; & t>0, \ 0\leq r<h(t),  \\
u_r(t,0)=0,\; u(t, h(t))=0,\quad & t>0, \\
h'(t)=-\mu u_r(t,h(t)),\quad &t>0,\\
h(0)=h_0,  u(0, r)=u_{0}(r),\; &0\leq r\leq h_0.
\end{array} \right.
\label{fa2}
\end{align}

\begin{prop}\label{DuLin1} Let $(u(t, r), h(t))$ be the solution of the  free boundary problem \eqref{fa2}.
Then the following alternative holds:

Either
\begin{itemize}
\item[(i)] {\rm Spreading:} $h_\infty=+\infty$ and $\lim_{t\to
+\infty}u(t, r)=\frac ab$ uniformly in any compact subset of $[0, \infty)$;
\end{itemize}

or
\begin{itemize}
\item[(ii)] {\rm Vanishing:} $h_\infty \leq R^* \sqrt{\frac da}$ and $\lim_{t\to +\infty}\|u(t,
\cdot)\|_{C([0,h(t)])}=0$.
\end{itemize}
\end{prop}

\begin{prop}\label{DuLin2} If $h_0 \geq R^*\sqrt{\frac da}$, then
spreading always happens.
If $h_0 < R^*\sqrt{\frac da}$, then there exists
$\mu^*>0$ depending on $u_0$ such that vanishing happens when $\mu\leq \mu^*$, and
spreading happens when $\mu>\mu^*$.
\end{prop}
The proofs of Propositions \ref{DuLin1} and \ref{DuLin2} can be found in \cite{DL} for the one dimensional case and \cite{DG} for higher
space dimensions.

\smallskip

We will show that when \eqref{u-sup} holds, a similar spreading-vanishing dichotomy holds for \eqref{f1}. More precisely, we have the following results.

\begin{thm}
\label{dichotomy}
Suppose that \eqref{u-sup} holds and $(u,v,h)$ is the unique solution of \eqref{f1} with $v_0\not\equiv 0$. Then the following alternative holds:

Either
\begin{itemize}
\item[(i)] {\rm Spreading of $u$:}\ \ $h_\infty=+\infty$ and $\lim_{t\to
+\infty}\big(u(t, r), v(t,r)\big)=\left(\frac{a_1}{b_1}, \, 0\right)$ uniformly in any compact subset of $[0, \infty)$;
\end{itemize}

or
\begin{itemize}
\item[(ii)] {\rm Vanishing of $u$:}\ \ $h_\infty \leq R^* \sqrt{\frac{d_1}{a_1-a_2c_1/c_2}}$\ and \ $\lim_{t\to +\infty}\|u(t,
\cdot)\|_{C([0,h(t)])}=0$, $\lim_{t\to\infty}v(t,r)=\frac{a_2}{c_2}$ uniformly in  any compact subset of $[0,\infty)$.
\end{itemize}
\end{thm}

\begin{thm}
\label{criteria}
In Theorem \ref{dichotomy}, if $h_0\geq  R^*\sqrt{\frac{d_1}{a_1-a_2c_1/c_2}}$, then spreading of $u$ always happens.
If $h_0< R^*\sqrt{\frac{d_1}{a_1-a_2c_1/c_2}}$, then there exists $\mu^*\in [0,\infty)$, depending on $(u_0,v_0)$, such that spreading of
$u$ happens exactly when $\mu>\mu^*$. Moreover, $\mu^*>0$ when $h_0< R^*\sqrt{\frac{d_1}{a_1}}$.
\end{thm}

We prove these results by several lemmas. In the rest of this section, we always assume that \eqref{u-sup} holds, and $(u,v,h)$ is the unique solution of
\eqref{f1}, with $v_0\not\equiv 0$.

Before starting the proofs, let us note that by the symmetric positions of $u$ and $v$ in \eqref{entire}, we may use Remark \ref{rmk-u-inf} to conclude that
when \eqref{u-sup} holds, the unique solution $(u,v)$ of the Cauchy problem of \eqref{entire} with $u(0,x), v(0,x)$ bounded, nonnegative and not identically zero, satisfies $\lim_{t\to\infty}\big(u(t,x), v(t,x)\big)=R_1=\left(\frac{a_1}{b_1},\, 0\right)$ locally uniformly in $\R^N$. The behavior of the free boundary problem
\eqref{f1} described by Theorems \ref{dichotomy} and \ref{criteria}, however, is more complicated, where vanishing for $u$ is possible.

\begin{lem}  If $h_\infty=\infty$, then $\lim_{t\to +\infty}u(t, r)=\frac
{a_1}{b_1}$ and  $\lim_{t\to +\infty}v(t, r)=0$ uniformly in any bounded subset of $[0, \infty)$.
\end{lem}
\bpf First we recall that the comparison principle gives
\begin{align*}
\limsup_{t\to +\infty}u(t, r)\leq \frac {a_1}{b_1}\ \ \textrm{uniformly for}\ r\in [0,\infty),\\
\limsup_{t\to +\infty}v(t, r)\leq \frac {a_2}{c_2}\ \ \textrm{uniformly for}\ r\in [0,\infty).
\end{align*}
Therefore for $\varepsilon_2=(\frac {a_1}{c_1}-\frac {a_2}{c_2})/2$, there exists $t_2>0$ such that $v(t, r)\leq \frac {a_2}{c_2}+\varepsilon_2$ for
$t\geq t_2, r\in [0,\infty)$.
 Therefore $u$ satisfies
  \begin{align}
\left\{
\begin{array}{lll}
u_{t}-d_1 \Delta u\geq u(c_1\varepsilon_2 -b_1 u),\; & t>t_2, \ 0\leq r<h(t),  \\
u_r(t,0)=0,\; u(t, h(t))=0,\quad & t>t_2,\\
h'(t)=-\mu u_r(t, h(t)),&t>t_2,\\
u(t_2, r)>0, &0\leq r<h(t),\\
\end{array} \right.
\label{fv2}
\end{align}

Next, we prove that for any $l> R^*\sqrt {\frac {d_1}{c_1\varepsilon_2}}$,
there exists $t_l>t_2$ such that
\begin{align}
\label{ine2}
u(t, r)\geq \frac{c_1\varepsilon_2}{2b_1}\  \textrm{ for}\ t\geq t_l,\ 0\leq r\leq l.
\end{align}
Indeed, since $h(t)\to \infty$, there exists $t_3\geq t_2$ such that $h(t_3)\geq l$. It follows from the
comparison principle (\cite{DG}) that $u(t, r)\geq \underline u_l(t,r)$ for $t>t_3$, $r\in [0, \underline h(t)]$ and $h(t)\geq \underline h(t)$ in
$(t_3, \infty)$, where $(\underline u_l, \underline h(t))$ is the solution
of the following free boundary problem discussed in \cite{DG}:
\begin{align} \left\{
\begin{array}{ll}
(\underline u_l)_t-d_1 \Delta \underline u_l=\underline u_l(c_1\varepsilon_2-b_1\underline u_l),\; &t>t_3,\ 0\leq r<\underline h(t),  \\
(\underline u_l)_{r}(t, 0)=\underline u_l(t, \underline h(t))=0, &t>t_3,\\
\underline h'(t)=-\mu\underline u_r(t, \underline h(t)),&t>t_3,\\
\underline u_l(t_3, r)=u(t_3, r),\; \underline h(t_3)=h(t_3), &0\leq r\leq h(t_3).
\end{array} \right.
\label{g2}
\end{align}
Since  $c_1\varepsilon_2> d_1(\frac {R^*}{l})^2$, it follows from Proposition \ref{DuLin2} that $\underline h(t)\to \infty$ and $\underline
u_l (t, r) \to \frac{c_1\varepsilon_2}{b_1}$ as $t\to \infty$ uniformly in  $[0, l]$, which implies that (\ref{ine2}) holds for some $t_l$.

Now we know that $(u, v)$ satisfies
\begin{align}
\left\{
\begin{array}{lll}
u_{t}-d_1 \Delta u=u(a_1-b_1 u-c_1 v),\; & t>t_l, \ 0\leq r<l,  \\
v_{t}-d_2 \Delta v=v(a_2-b_2 u-c_2 v),\; & t>t_l, \ 0\leq r<l,  \\
u_r(t,0)=v_r(t,0)=0,\quad & t>t_l,\\
u(t, r)\geq \frac{c_1\varepsilon_2}{2b_1},\ v(t, r)\leq \frac {a_2}{c_2}+\varepsilon_2,  & t\geq t_l,\ 0\leq r\leq l.
\end{array} \right.
\label{f2s}
\end{align}
As in the proof of Theorem \ref{vanofu}, it follows from the theory of monotone dynamical systems  that
$\liminf_{t\to +\infty}u(t, r)\geq \underline u_l(r)$ and $\limsup_{t\to +\infty}v(t, r)\leq \overline v_l(r)$ in
$[0, l]$, where $(\underline u_l, \overline v_l)$  satisfies
\begin{align*} \left\{
\begin{array}{lll}
-d_1 \Delta \underline u_l=\underline u_l(a_1-b_1 \underline u_l-c_1\overline v_l),\; &\ 0\leq r<l,  \\
-d_2 \Delta \overline v_l=\overline v_l(a_2-b_2 \underline u_l-c_2\overline v_l),\; &\ 0\leq r<l, \\
\frac{\partial \underline u_l}{\partial r}(0)=\frac{\partial \overline v_l}{\partial r}(0)=0,\; & \\
\underline u_l(l)=\frac{c_1\varepsilon_2}{2b_1},\ \overline v_l(l)=\frac {a_2}{c_2}+\varepsilon_2. &
\end{array} \right.
\end{align*}
Letting $l\to \infty$, we similarly have  $(\underline u_l, \overline v_l)\to(\underline u_\infty, \overline v_\infty)$ which satisfies
\begin{align*} \left\{
\begin{array}{lll}
-d_1 \Delta \underline u_\infty=\underline u_\infty(a_1-b_1 \underline u_\infty-c_1\overline v_\infty),\; &\ 0\leq r<\infty,  \\
-d_2 \Delta \overline v_\infty=\overline v_\infty(a_2-b_2 \underline u_\infty-c_2\overline v_\infty),\; &\ 0\leq r<\infty, \\
\frac{\partial \underline u_\infty}{\partial r}(0)=\frac{\partial \overline v_\infty}{\partial r}(0)=0,\; & \\
\underline u_\infty(r)\geq \frac{c_1\varepsilon_2}{2b_1},\ \overline v_\infty(r)\leq \frac {a_2}{c_2}+\varepsilon_2, &0\leq r<\infty.
\end{array} \right.
\end{align*}

 Due to \eqref{u-sup}, as before, by  the global dynamical behavior of the associated  ODE system (\cite{MT}), we deduce
  that $\underline u_\infty(r)=\frac{a_1}{b_1}$ and $\overline v_\infty(r)=0$.
We thus obtain $\liminf_{t\to +\infty}u(t, r)\geq \frac {a_1}{b_1}$ and $\limsup_{t\to +\infty}v(t, r)\leq 0$
uniform in $[0,l]$, which implies that
 $\lim_{t\to +\infty}u(t, r)=\frac {a_1}{b_1}$
and  $\lim_{t\to +\infty}v(t, r)=0$ uniformly in any bounded subset of $[0, \infty)$.
\epf

\begin{lem}   If $h_\infty<\infty$, then $\lim_{t\to
+\infty}\|u(t, \cdot)\|_{C([0, h(t)])}=0$ and  $\lim_{t\to +\infty}v(t, r)=\frac
{a_2}{c_2}$ uniformly in any bounded subset of $[0, \infty)$.
\end{lem}
\bpf
Define
\begin{align*}
s=\frac{h_0r}{h(t)},\ w(t, s)=u(t, r), \ z(t, s)=v(t, r).
\end{align*}
Direct calculations yield
\begin{align*}
u_t=w_t-\frac{h'(t)}{h(t)}sw_s,\ u_r=\frac{h_0}{h(t)}w_s,\ \Delta_r u=\frac{h^2_0}{h^2(t)}\Delta_s w.
\end{align*}
Hence $w(s, t)$ satisfies
\begin{align}
\left\{
\begin{array}{lll}
w_{t}-d_1\frac{h_0^2}{h^2(t)}\Delta_s w-\frac{h'(t)}{h(t)}sw_s=w(a_1-b_1w-c_1z),\; &
0\leq s<h_0,\; t>0, \\
w_s(t, 0)=w(t, h_0)=0,&t>0,\\
w(0,s)=u_0(s)\geq 0,\; & 0\leq s\leq h_0.
\end{array} \right.
\label{Bb1}
\end{align}
This is an initial boundary value problem over a fixed ball $\{s<h_0\}$. Since $h_0\leq h(t)<h_\infty<\infty$,
the differential operator is uniformly parabolic. By Theorem 2.5,
\begin{align*}
\|a_1-b_1w-c_1z\|_{L^\infty}\leq a_1+b_1M_1+c_2M_2,\ \left\|\frac{h'(t)}{h(t)}s\right\|_{L^\infty}\leq M_3.
\end{align*}
Therefore we can apply standard $L^p$ theory  and then the Sobolev imbedding
theorem (\cite{LSU}) to obtain, for any $\alpha\in (0,1)$,
\begin{align*}
\|w\|_{C^{(1+\alpha)/2,1+\alpha}([0,\infty)\times[0,
h_0])}\leq  C_4,
\end{align*}
where $C_4$ is a constant depending on
$\alpha, h_0, M_1, M_2, M_3$ and $\|u_0\|_{C^{1+\alpha}[0, h_0]}$. Similarly we may use interior estimates to the equation of $z$
to obtain
\begin{align*}
\|z\|_{C^{(1+\alpha)/2,1+\alpha}([0,\infty)\times[0,
h_0])}\leq  C_5,
\end{align*}
where $C_5$ is a constant depending on
$\alpha, h_0, M_1, M_2, M_3$ and $\|v_0\|_{C^{1+\alpha}[0, h_0+1]}$.

It follows that
there exists a constant $\tilde C$ depending on $\alpha, h_0,
(u_0,v_0)$ and $h_\infty$ such that
\begin{align}\|u\|_{C^{
(1+\alpha)/2,
1+\alpha}(G)}+\|v\|_{C^{
(1+\alpha)/2,
1+\alpha}(G)}+\|h\|_{C^{1+\alpha/2}([0,\infty))}\leq
\tilde C,\label{Bg1}
\end{align}
where $G:=\{(t,r): t\geq 0, r\in [0, h(t)]\}$.

Arguing indirectly, we assume that  $\limsup_{t\to
+\infty}\|u(t, \cdot)\|_{C([0, h(t)])}=\delta>0$. Then there exists a sequence $(t_k, r_k)$ with $0<t_k<\infty$, $0\le r_k<h(t_k)$ such that $u(t_k, r_k)\geq \delta /2$ for all $k \in \mathbb{N}$, and $t_k\to \infty$ as $k\to \infty$.
Since $u(t, h(t))=0$ and since \eqref{Bg1} infers that $|u_r(t, h(t))|$ is uniformly bounded for $t\in [0,\infty)$,
there exists $\sigma>0$ such that $r_k\leq h(t_k)-\sigma$ for all $k\geq 1$. Therefore a subsequence of $\{r_k\}$ converges to $r_0\in [0, h_\infty-\sigma]$. Without loss of generality,
we assume $r_k\to r_0$ as $k\to \infty$.

Define
\[\mbox{
$u_k(t, r)=u(t_k+t,r)$ and $v_k(t, r)=v(t_k+t, r)$ for $(t, r)\in G_k$}
\]
with
\[
G_k:=\{(t,r): t\in(-t_k, \infty), r\in [0, h(t_k+t)]\}.
\]
It follows from \eqref{Bg1} that  $\{(u_k, v_k)\}$ has a subsequence $\{(u_{k_i}, v_{k_i})\}$ such that
\[
\mbox{$\|(u_{k_i}, v_{k_i})-(\tilde u, \tilde v)\|_{C^{1,2}(G_{k_i})\times C^{1,2}(G_{k_i})}\to 0$ as $i\to \infty$,}
\] and $(\tilde u, \tilde v)$ satisfies
\begin{align*} \left\{
\begin{array}{lll}
\tilde u_t-d_1 \Delta \tilde u=\tilde u(a_1-b_1 \tilde u-c_1\tilde v),\; &\ t\in (-\infty, \infty),\ 0\leq r<h_\infty,  \\
\tilde v_t-d_2 \Delta \tilde v=\tilde v(a_2-b_2 \tilde u-c_2\tilde v),\; &\ t\in (-\infty, \infty),\ 0\leq r<h_\infty,\\
\end{array} \right.
\end{align*}
with
\[
\tilde u(t,h_\infty)=0 \;\mbox{ for } t\in (-\infty,\infty).
\]
Since $\tilde u(0, r_0)\geq \delta/2$, the maximum principle infers that $\tilde u>0$ in $(-\infty, \infty)\times [0, h_\infty)$.
Thus we can apply the Hopf boundary lemma to conclude
that $\sigma_0:=\tilde u_r(0, h_\infty)<0$. It follows that
\[
u_r(t_{k_i}, h(t_{k_i}))=\partial_ru_{k_i}(0, h(t_{k_i}))\leq \sigma_0/2<0
\]
for all large $i$, and hence $h'(t_{k_i})\geq -\mu\sigma_0/2>0$ for all large $i$.

On the other hand,
since $\|h\|_{C^{1+\alpha/2}([0,\infty))}\leq
\tilde C$, $h'(t)>0$ and $h(t)\leq h_\infty$,  we necessarily  have $h'(t)\to 0$ as $t\to \infty$.
 This contradiction shows that we must have
 \[
 \lim_{t\to +\infty}\|u(t,
\cdot)\|_{C([0,h(t)])}=0.
\]
 We may now use a simple comparison argument to deduce that $\lim_{t\to +\infty}v(t, r)=\frac
{a_2}{c_2}$ uniformly in any bounded subset of $[0, \infty)$.
\epf

\begin{lem}\label{vanofu3}
 If $h_\infty<\infty$, then $h_\infty\leq R^*\sqrt{\frac {d_1}{a_1-a_2c_1/c_2}}$. Hence $h_0\geq R^*\sqrt{\frac {d_1}{a_1-a_2c_1/c_2}}$
 implies  $h_\infty=\infty$.
\end{lem}
\bpf  Assume for contradiction that  $ R^*\sqrt{\frac {d_1}{a_1-a_2c_1/c_2}}<h_\infty<\infty$. Then there exists
$T_1>0$ such that $l:=h(T_1)>R^*\sqrt{\frac {d_1}{a_1-a_2c_1/c_2}}$. We choose $\varepsilon$ sufficiently small such that
$ h(T_1)>R^*\sqrt{\frac {d_1}{a_1-a_2c_1/c_2-c_1\varepsilon}}$.

Recall that
\begin{align*}
\limsup_{t\to +\infty}v(t, r)\leq \frac {a_2}{c_2}\ \ \textrm{uniformly for}\ r\in [0,\infty).
\end{align*}
Therefore for the above chosen $\varepsilon>0$, there exists $T_2>T_1$ such that $v(t, r)\leq \frac {a_2}{c_2}+\varepsilon$ for
$t\geq T_2, r\in [0,\infty)$.
 Hence $(u, h)$ satisfies
  \begin{align}
\left\{
\begin{array}{ll}
u_{t}-d_1 \Delta u\geq u(a_1-c_1(a_2/c_2+\varepsilon) -b_1 u),\; & t>T_2, \ 0\leq r<h(t),  \\
u_r(t,0)=0,\; u(t, h(t))=0,\quad & t>T_2,\\
h'(t)=-\mu u_r(t, h(t)),&t>T_2,\\
u(T_2, r)>0, &0\leq r<h(T_2).
\end{array} \right.
\label{fg1}
\end{align}
This indicates that $(u, h)$ is an upper solution to the problem
\begin{align}
\left\{
\begin{array}{ll}
\underline u_{t}-d_1  \Delta \underline u=\underline u(a_1-c_1(a_2/c_2+\varepsilon) -b_1 \underline u),\; & t>T_2, \ 0\leq r<\underline h(t),  \\
\underline u_r(t,0)=0,\; \underline u(t, \underline h(t))=0,\quad & t>T_2,\\
\underline h'(t)=-\mu \underline u_r(t, \underline h(t)),&t>T_2,\\
\underline u(T_2, r)=u(T_2, r),\  \underline h(T_2)=h(T_2), &0\leq r<\underline h(T_2).
\end{array} \right.
\label{fg2}
\end{align}
Hence $h(t)\geq \underline h(t)$ for $t>T_2$.
Since $\underline h(T_2)> l>R^*\sqrt{\frac {d_1}{a_1-a_2c_1/c_2-c_1\varepsilon}}$, it follows from Proposition 4.2 that
$\underline h(t)\to \infty$ and therefore $h_\infty=\infty$.
 This contradiction proves that $h_\infty\leq R^*\sqrt{\frac {d_1}{a_1-a_2c_1/c_2}}$.
 \epf

\begin{lem}\label{spread2} If $h_0<R^*\sqrt{\frac {d_1}{a_1-a_2c_1/c_2}}$, then there exists
$\underline \mu\geq 0$ depending on $u_0$ and $v_0$ such that
$h_\infty=\infty$ when $\mu>\underline \mu$.
\end{lem}
\bpf  Since
\begin{align*}
\limsup_{t\to +\infty}v(t, r)\leq \frac {a_2}{c_2}\ \textrm{uniformly for}\ r\in [0,\infty),
\end{align*}
 for $\varepsilon=(\frac {a_1}{c_1}-\frac {a_2}{c_2})/2>0$, there exists $T$, which is independent of $\mu$,
 such that $v(t, r)\leq \frac {a_2}{c_2}+\varepsilon$ for
$t\geq T, r\in [0,\infty)$.
 Thus $(u, h)$ satisfies
  \begin{align}
\left\{
\begin{array}{lll}
u_{t}-d_1 \Delta u\geq u(c_1\varepsilon -b_1 u),\; & t>T, \ 0\leq r<h(t),  \\
u_r(t,0)=0,\; u(t, h(t))=0,\quad & t>T,\\
h'(t)=-\mu u_r(t, h(t)),&t>T,\\
u(T, r)>0, &0\leq r<h(T).
\end{array} \right.
\label{fu3}
\end{align}
Note that $u(T,r)$ depends on $\mu$, but $u(T,r)\geq z(T,r)$ in $[0,\infty)\times [0, h_0]$, where $z(t,r)$ and $w(t,r)$ are determined by
\begin{align*}\left\{
\begin{array}{ll}
z_t-d_1 \Delta z=z(a_1-b_1 z-c_1 w),& t>0, \ 0\leq r<h_0,  \\
w_t-d_2 \Delta w=w(a_2-b_2 z-c_2 w),& t>0, \ 0\leq r<h_0, \\
z_r(t,0)=w_r(t,0)=0,\;  & t>0,\\
z(t, h_0)=0,\; w(t, h_0)=\max\{\frac{a_2}{c_2},\|v_0\|_{L^{\infty}}\},& t>0,\\
z(0, r)=u_{0}(r),\ w(0, r)=\max\{\frac{a_2}{c_2},\|v_0\|_{L^{\infty}}\}, &0\leq r\leq h_0.
\end{array} \right.
\end{align*}
Clearly, $z(T,r)$ is independent of $\mu$. Now it is easy to see that $(u, h)$ is an upper solution to
the problem
\begin{align}
\left\{
\begin{array}{lll}
\underline u_{t}-d_1  \Delta \underline u=\underline u(c_1\varepsilon -b_1 \underline u),\; & t>T, \ 0\leq r<\underline h(t),  \\
\underline u_r(t,0)=0,\; \underline u(t, \underline h(t))=0,\quad & t>T,\\
\underline h'(t)=-\mu \underline u_r(t, \underline h(t)),&t>T,\\
\underline u(T, r)=z(T, r),\  \underline h(T)=h_0, &0\leq r\leq h_0.
\end{array} \right.
\label{ff1}
\end{align}
By Lemma 2.8 of \cite{DG}, there exists $\underline\mu>0$ such that $\underline h_\infty=+\infty$ for $\mu>\underline\mu$.
One can actually  argue as in Lemma 3.7 of \cite{DL} to show that if
\begin{align*}
\mu \geq
\underline \mu:=\max\left\{1,\frac{b_1}{c_1\varepsilon}\|z(T,r)\|_\infty\right\}\frac {d_1((R^*\sqrt{\frac {d_1}{a_1-a_2c_1/c_2}})^N-h^N_0)}{N(\int ^{h_0}_0 r^{N-1}z(T,r)dr)},
\end{align*}
then $\underline h_\infty= +\infty$. Therefore $h_\infty=\infty$ for $\mu >\underline \mu$.
 \epf

\begin{lem}  There exists $\mu^*\geq 0$ depending on $u_0$ and $v_0$ such that
$h_\infty =+\infty$ if $\mu>\mu^*$ and $h_\infty <\infty$ if $0<\mu\leq \mu^*$.
\end{lem}
\bpf Define $\Sigma :=\left\{\mu>0: h_\infty > R^*\sqrt{\frac {d_1}{a_1-a_2c_1/c_2}}\right\}$. It follows from Lemmas \ref{vanofu3} and \ref{spread2} that
 $\mu^*:=\inf
\Sigma\in [0, \infty)$. By Lemma \ref{vanofu3} and the monotonicity of $h_\infty$ with respect to $\mu$ (Corollary \ref{monotone}),
we find that $h_\infty=+\infty$ when $\mu>\mu^*$ and $h_\infty<+\infty$ when $0<\mu<\mu^*$.

We claim that if $\mu^*>0$, then $\mu^*\not\in\Sigma$. Otherwise $h_\infty > R^*\sqrt{\frac {d_1}{a_1-a_2c_1/c_2}}$ for
$\mu=\mu^*$. Hence we can find $T>0$ such that $h(T)> R^*\sqrt{\frac {d_1}{a_1-a_2c_1/c_2}}$. To stress the dependence of the solution $(u,v,h)$
of \eqref{f1} on $\mu$, we now write $(u_\mu,v_\mu, h_\mu)$ instead of
$(u,v,h)$. So we have $h_{\mu^*}(T)> R^*\sqrt{\frac {d_1}{a_1-a_2c_1/c_2}}$. By
the continuous dependence of $(u_\mu, v_\mu, h_\mu)$ on $\mu$, we can find
$\epsilon>0$ small so that $h_\mu(T)> R^*\sqrt{\frac {d_1}{a_1-a_2c_1/c_2}}$ for
all $\mu\in [\mu^*-\epsilon, \mu^*+\epsilon]$. It follows that for
all such $\mu$,
\begin{align*}
\lim_{t\to\infty} h_\mu(t)>h_\mu(T)>R^*\sqrt{\frac {d_1}{a_1-a_2c_1/c_2}} \, .
\end{align*}
This implies that $[\mu^*-\epsilon, \mu^*+\epsilon]\subset
\Sigma$, and $\inf\Sigma\leq \mu^*-\epsilon$,
contradicting the definition of $\mu^*$. This proves our claim that
$\mu^*\in\Sigma$.
 \epf

Lemma \ref{vanofu3} implies that the threshold constant $\mu ^*=0$ if $h_0\geq R^*\sqrt{\frac {d_1}{a_1-a_2c_1/c_2}}$.
The next lemma shows that $\mu ^*>0$ if $h_0< R^*\sqrt{\frac
{d_1}{a_1}}$.

\begin{lem} Suppose $h_0< R^*\sqrt{\frac
{d_1}{a_1}}$.  Then there exists $\overline \mu>0$ depending on $u_0$ such
that $h_\infty <+\infty$ if $\mu\leq \overline \mu$.
\end{lem}
\bpf Clearly
$(u, h)$ satisfies
  \begin{align}
\left\{
\begin{array}{lll}
u_{t}-d_1 \Delta u\leq u(a_1-b_1 u),\; & t>0, \ 0\leq r<h(t),  \\
u_r(t,0)=0,\; u(t, h(t))=0,\quad & t>0,\\
h'(t)=-\mu u_r(t, h(t)),&t>0,\\
u(0, r)=u_0(r), &0\leq r<h_0.
\end{array} \right.
\label{fg3}
\end{align}
This indicates that $(u, h)$ is a lower solution to the problem
\begin{align}
\left\{
\begin{array}{lll}
\overline u_{t}-d_1  \Delta \overline u=\overline u(a_1-b_1 \overline u),\; & t>0, \ 0\leq r<\overline h(t),  \\
\overline u_r(t,0)=0,\; \overline u(t, \overline h(t))=0,\quad & t>0,\\
\overline h'(t)=-\mu \overline u_r(t, \overline h(t)),&t>0,\\
\overline u(0, r)=u_0(r),\  \overline h(0)=h_0, &0\leq r<\overline h(0).
\end{array} \right.
\label{f1g2}
\end{align}
Since  $h_0< R^*\sqrt{\frac
{d_1}{a_1}}$, it follows from Proposition \ref{DuLin2} that there exists $\overline \mu>0$ depending on $u_0$ such
that $\overline h_\infty <+\infty$ if $\mu\leq \overline \mu$, therefore $h_\infty <+\infty$ if $\mu\leq \overline \mu$.
 \epf

\smallskip

Theorems \ref{dichotomy} and \ref{criteria} now follow directly from the conclusions proved in the above lemmas.

\section{estimates of spreading speed}
In this section we give some rough estimates on the spreading speed of $h(t)$ for the case that spreading of $u$ happens.
We always assume that \eqref{u-sup} holds.

We first recall Proposition 3.1 of \cite{DG}, whose complete proof is given in \cite{BDK}.

\begin{prop} For any given constants $a>0$,
$b>0$, $d>0$ and $k \in [0,2\sqrt{ad})$, the problem
\begin{equation}
\label{4.1} -dU''+k U'=a U-b U^2 \quad \mbox{ in } (0,\infty),
\;\;\; U(0)=0
\end{equation}
admits a unique positive solution $U=U_k=U_{a,b,d,k}$, and it
satisfies $U(r) \to \frac{a}{b}$ as $r \to +\infty$. Moreover,
$U'_{k} (r)> 0$ for $r \geq 0$, $U_{k_1}'(0)>U_{k_2}'(0)$,
$U_{k_1}(r)> U_{k_2}(r)$ for $r>0$ and $k_1< k_2$, and for each
$\mu>0$, there exists a unique $k_0=k_0(\mu,a,b, d)\in (0, 2\sqrt{ad})$
such that $\mu U'_{k_0}(0)=k_0$. Furthermore,
\begin{equation}
\label{k0-lim} \lim_{\frac{a\mu}{bd}\to
\infty}\frac{k_0}{\sqrt{ad}}=2,\; \lim_{\frac{a\mu}{bd}\to
0}\frac{k_0}{\sqrt{ad}}\frac{bd}{a\mu}=1/\sqrt 3.
\end{equation}
\end{prop}

Making use of the function $k_0(\mu,a,b,d)$, we have the following estimates for the spreading speed of $h(t)$.

\begin{thm}\label{spreadsp1} Assume that \eqref{u-sup} holds.
 If $h_\infty=+\infty$, then
\begin{align*}
k_0(\mu, a_1-a_2c_1/c_2, b_1, d_1)\leq \liminf_{t\to +\infty} \frac {h(t)}t\leq \limsup_{t\to +\infty} \frac {h(t)}t\leq k_0(\mu, a_1, b_1, d_1).
\end{align*}
\label{m10}
\end{thm}
\bpf  Since
\[
\left\{
\begin{array}{ll}
u_t-d_1 \Delta u=u(a_1-b_1u-c_1v)\leq u(a_1-b_1u), & t>0, 0\leq r<h(t),\\
u_r(t,0)=0, \; u(t, h(t))=0, & t>0,\\
h'(t)=-\mu u_r(t, h(t)), & t>0,\\
u(0,r)=u_0(r),& 0\leq r\leq h_0,
\end{array}\right.
\]
the pair $(u,h)$ is a lower solution to the problem
\[\left\{
\begin{array}{ll}
\overline u_t-d_1 \Delta \overline u=\overline  u(a_1-b_1\overline u), & t>0, 0\leq r<\overline h(t),\\
\overline u_r(t,0)=0, \; \overline u(t, \overline h(t))=0, & t>0,\\
\overline h'(t)=-\mu \overline u_r(t, \overline h(t)), & t>0,\\
\overline u(0,r)=u_0(r),& 0\leq r\leq h_0,
\end{array}
\right.
\]
It follows that $\overline h(t)\geq h(t)\to\infty$ as $t\to\infty$. By \cite{DG},
\[
\lim_{t\to\infty} \frac{\overline h(t)}{t}=k_0(\mu, a_1, b_1, d_1).
\]
We thus have
\[
\limsup_{t\to +\infty} \frac {h(t)}t\leq k_0(\mu, a_1, b_1, d_1).
\]

Next we prove that
\[
k_0(\mu,a_1-a_2c_1/c_2, b_1, d_1)\leq \liminf_{t\to +\infty} \frac {h(t)}t.
\]
Since $\limsup_{t\to\infty}v\leq \frac {a_2}{c_2}$ uniformly for
$r\in [0,\infty)$ and $h_\infty=\infty$,  for any $0<\varepsilon <\varepsilon_0:=(\frac {a_1}{c_1}-\frac {a_2}{c_2})/2$,
there exists $T_\varepsilon$, such that $v(t, r)\leq \frac {a_2}{c_2}+\varepsilon$  for
$t\geq T_\varepsilon, r\in [0,\infty)$ and $h(T_\varepsilon)>R^*\sqrt{\frac {d_1}{a_1-c_1(a_2/c_2+\varepsilon)}}$.
 Hence $(u, h)$ satisfies
  \begin{align}
\left\{
\begin{array}{lll}
u_{t}-d_1 \Delta u\geq u(a_1-c_1(a_2/c_2+\varepsilon)-b_1 u),\; & t>T_\varepsilon, \ 0\leq r<h(t),  \\
u_r(t,0)=0,\; u(t, h(t))=0,\quad & t>T_\varepsilon,\\
h'(t)=-\mu u_r(t, h(t)),&t>T_\varepsilon,\\
u(T_\varepsilon, r)>0, &0\leq r<h(T_\varepsilon).
\end{array} \right.
\label{f51}
\end{align}
This implies that $(u, h)$ is an upper solution to the problem
\begin{align}
\left\{
\begin{array}{lll}
\underline u_{t}-d_1  \Delta \underline u=\underline u(a_1-c_1(a_2/c_2+\varepsilon) -b_1 \underline u),\; & t>T_\varepsilon, \ 0\leq r<\underline h(t),  \\
\underline u_r(t,0)=0,\; \underline u(t, \underline h(t))=0,\quad & t>T_\varepsilon,\\
\underline h'(t)=-\mu \underline u_r(t, \underline h(t)),&t>T_\varepsilon,\\
\underline u(T_\varepsilon, r)=u(T_\varepsilon, r), &0\leq r\leq h(T_\epsilon).
\end{array} \right.
\label{f52}
\end{align}
By Proposition \ref{DuLin1}, $\underline h_\infty=\infty$ since  $h(T_\varepsilon)>R^*\sqrt{\frac {d_1}{a_1-c_1(a_2/c_2+\varepsilon)}}$.
  Moreover, from \cite{DG} we have $\lim_{t\to \infty}\frac {\underline h(t)}{t}=k_0(\mu, a_1-c_1(a_2/c_2+\varepsilon), b_1, d_1)$,
which implies that
\[\mbox{
 $k_0(\mu, a_1-c_1(a_2/c_2+\varepsilon), b_1, d_1)\leq \liminf_{t\to +\infty} \frac {h(t)}t$ for
any $0<\varepsilon <\varepsilon_0$.}
\]
 Letting $\varepsilon \to 0$ and using the continuity of $k_0$ with respect to its arguments,
we immediately obtain  the  desired result.
\epf

Next we obtain an upper bound for the spreading speed that is independent of $\mu$, namely, we show that under suitable additional conditions on $(u_0,v_0)$,
\begin{equation}
\label{speed-ub}
\limsup_{t\to\infty}\frac{h(t)}{t}\leq c_*,
\end{equation}
where $c_*$ is the minimal speed of the traveling waves to \eqref{entire} in dimension one given in \cite{KO}.
More precisely, by Theorem 2.1 of \cite{KO}, when \eqref{u-sup} holds, there exists $c_*>0$, depending continuously  on
the parameters in \eqref{entire}, such that \eqref{entire} with $N=1$ has a solution of the form \eqref{tw} when $c\geq c_*$, and there
is no such solution when $c<c_*$. (We remark that by a standard change of variables, our general form \eqref{entire} can be reduced to
the special form considered in \cite{KO}.)

\begin{thm}
\label{thm-speed-ub}
Suppose that \eqref{u-sup} holds, $h_\infty=\infty$, and additionally
\begin{equation}
\label{initial}\mbox{
$u_0(r)\leq \frac{a_1}{b_1}$ in $[0, h_0]$,\; $v_0(r)>0$ in $[0,\infty)$,\;
 $\liminf_{r\to\infty}v_0(r)\geq \frac{a_2}{c_2}$.
 }
 \end{equation}
  Then \eqref{speed-ub} holds.
\end{thm}

\bpf
In \eqref{entire}, we replace $a_1$ by  $\tilde a_1=a_1+\epsilon_1$ and replace $a_2$ by $\tilde a_2=a_2-\epsilon_2$,
where $\epsilon_1,\; \epsilon_2$ are  small positive constants. We then denote by $\tilde c_*$ the minimal speed of the traveling waves
to the modified \eqref{entire}, and by $(U(x),V(x))$ the corresponding solution in \eqref{tw}. We thus have
\[
\left\{
\begin{array}{ll}
U''+\tilde c_*U'+U(\tilde a_1-b_1U-c_1V)=0,\;\; U'<0,& x\in \R^1,\\
V''+\tilde c_* V'+V(\tilde a_2-b_2U-c_2V)=0,\;\; V'>0, & x\in\R^1,\\
U(-\infty)=\frac{\tilde a_1}{b_1},\; V(-\infty)=0,\; U(+\infty)=0,\; V(+\infty)=\frac{ \tilde a_2}{c_2}.&
\end{array}
\right.
\]

Let $\sigma(x)$ be a smooth function satisfying
\[
\sigma(x)\in [0,1],\; \sigma(x)=0 \mbox{ for } x\leq -1,\; \sigma(x)=1 \mbox{ for } x\geq 0,\; \sigma'(x)\geq 0.
\]
For $x_0>0$ and $L>0$ to be determined later, we define
\[
\tilde U(x)=\left\{\begin{array}{ll}
U(x)-\sigma\big(\frac{x-x_0+1}{L}\big)U(x_0), & x\leq x_0,\\
0, & x>x_0.
\end{array}
\right.
\]
Clearly $\tilde U(x)\leq U(x)$ and
\[
\left\{
\begin{array}{ll}
 \tilde U'(x)\leq 0, &x\in\R^1,\\
 \tilde U(x)=U(x)-U(x_0), &x\in [x_0-1, x_0],\\
\tilde U(x)=U(x), & x\leq x_0-1-L.
\end{array}
\right.
\]

We now define
\[
\xi(t) =\tilde c_*t+\xi_0,\;
u^*(t,r)=\tilde U(r-\xi(t)+x_0),\; v_*(t,r)=V(r-\xi(t)+x_0).
\]
We will show that by choosing $x_0,L, \xi_0$  properly, we have $h(t)\leq \xi(t)$ for all $t>0$.
 Clearly this implies
 \[
 \limsup_{t\to\infty}\frac{h(t)}{t}\leq \lim_{t\to\infty}\frac{\xi(t)}{t}=\tilde c_*.
 \]
 Letting  $\epsilon_1,\; \epsilon_2\to 0$, due to the continuous dependence of $\tilde c_*$ on the parameters, we
 deduce
 \[
 \limsup_{t\to\infty}\frac{h(t)}{t}\leq c_*,
 \]
 as desired.

 It remains to show $h(t)\leq \xi(t)$ for $t>0$. This will be shown by checking that
 $(u^*,v_*,\xi)$ can serve as an upper solution to \eqref{f1}. Clearly
 \[
 u^*(t,\xi(t))=\tilde U(x_0)=0 \;\;\forall t>0.
 \]
 Moreover,
 \[-\mu u^*_r(t,\xi(t))=\mu \tilde U'(x_0)=\mu U'(x_0)<\tilde c_*=\xi'(t)
 \]
 provided $x_0>0$ is chosen large enough, since $U'(x)\to 0$ as $x\to+\infty$.

By direct calculations we obtain
\begin{align*}
&u^*_t-d_1\Delta u^*\\
&=-\tilde c_*\tilde U'(r-\xi(t)+x_0)-d_1\left(\tilde U''(r-\xi(t)+x_0)+\frac{N-1}{r}\tilde U'(r-\xi(t)+x_0)\right)\\
&\geq -\tilde c_*\tilde U'-d_1\tilde U''\\
&=-\tilde c_* U'-d_1 U''+U(x_0)\left(\tilde c_* \frac{\sigma'}{L}+d_1\frac{\sigma''}{L^2}\right)\\
&=U(\tilde a_1-b_1U-c_1V)+L^{-1}U(x_0)\left(\tilde c_* \sigma'+d_1\frac{\sigma''}{L}\right)\\
&=U(a_1-b_1U-c_1V)+\epsilon_1 U+L^{-1}U(x_0)\left(\tilde c_* \sigma'+d_1\frac{\sigma''}{L}\right).
\end{align*}

For $r\in (0, \xi(t)]$, we have $r-\xi(t)+x_0\leq x_0$ and therefore
$U(r-\xi(t)+x_0)\geq U(x_0)$. On the other hand, by choosing $L$ suitably large we can guarantee that
$L^{-1}\left|\tilde c_* \sigma'+d_1\frac{\sigma''}{L}\right|\leq \epsilon_1$. Therefore, with $L$ chosen this way, we have
\[
\epsilon_1 U+L^{-1}U(x_0)\left(\tilde c_* \sigma'+d_1\frac{\sigma''}{L}\right)\geq 0 \mbox{ for } r\in [0, \xi(t)].
\]
Moreover, we notice that $\tilde U(x)=U(x)$ for $x\leq x_0-1-L$, and by enlarging $x_0$ if necessary, we can make $(U(x_0-1-L),V(x_0-1-L))$ very close to
$(0, \frac{\tilde a_2}{c_2})$. It follows that for $x\in [x_0-1-L, x_0]$, $(\tilde U(x), V(x))$ is very close to $(0, \frac{\tilde a_2}{c_2})$.
Hence, due to \eqref{u-sup}, for every fixed $x\in [x_0-1-L, x_0]$, the function
\[
\eta(s):=s[a_1-b_1s-c_1V(x)]
\mbox{ is increasing in $[0, U(x)]$.}
\]
We may now use $U(x)\geq \tilde U(x)$ to deduce that, for $x\in [x_0-1-L, x_0]$,
\[
U(x)[a_1-b_1U(x)-c_1V(x)]\geq \tilde U(x)[a_1-b_1\tilde U(x)-c_1V(x)].
\]
We thus  obtain, by combining the above inequalities,
\[
u^*_t-d_1\Delta u^*\geq \tilde U(a_1-b_1\tilde U-c_1V)= u^*(a_1 -b_1u^*-c_1v_*)
\]
for $r\in (0,\xi(t)]$ and $t>0$.

We also have
\begin{align*}
(v_*)_t-d_2\Delta v_*&=-\tilde c_* V'-d_2 V''-d_2\frac{N-1}{r}V'\\
&\leq -\tilde c_* V'-d_2 V''\\
&=V(\tilde a_2-b_2U-c_2V)\\
&\leq V(\tilde a_2-b_2\tilde U-c_2V)\\
&\leq v_*(a_2-b_2u^*-c_2v_*)
\end{align*}
for $r>0$ and $t>0$.

For $r\in [0, h_0]$, $u^*(0,r)=\tilde U(r-\xi_0+x_0)\geq \tilde U(h_0-\xi_0+x_0)$.
Since $\tilde U(x)\to \frac{\tilde a_1}{b_1}$ as $x\to-\infty$, we may choose $\xi_0>h_0$ large enough such that
$\tilde U(h_0-\xi_0+x_0)\geq \frac{ a_1}{b_1}$ and hence, with such a $\xi_0$,
\[
u^*(0,r)\geq u(0,r),\; \forall r\in [0,h_0].
\]

By our assumption on
$v_0(r)$, there exists $X_1>0$ such that $v_0(r)\geq \frac{\tilde a_2}{c_2}$ for $r\geq X_1$.
Let $\sigma_0:=\inf_{r\in [0, X_1]}v_0(r)$. Then $\sigma_0>0$ and we can find $X_2>0$ such that
$V(x)<\sigma_0$ for $x\leq -X_2$. Therefore, for $r\in [0, X_1]$,
\[
v_*(0,r)=V(r-\xi_0+x_0)\leq V(X_1-\xi_0+x_0)<\sigma_0\leq v_0(r)
\]
if we enlarge $\xi_0$ further to guarantee $X_1-\xi_0+x_0\leq -X_2$.
Since $v_*(0,r)=V(r-\xi_0+x_0)<\frac{\tilde a_2}{c_2}$ for every $r$, and $v_0(r)\geq \frac{\tilde a_2}{c_2}$ for $r\geq X_1$, we thus have
\[
v_*(0,r)\leq v_0(r) \;\;\forall r\geq 0.
\]

We are now in a position to apply Lemma \ref{comparison} to deduce $h(t)\leq \xi(t)$, $u(t,r)\leq u^*(t,r)$ and $v(t,r)\geq v_*(t,r)$,
except that the conditions $\partial_r u^*(t,0)=\partial_r v_*(t,0)=0$ are not satisfied. Instead, we have
\[
\label{r=0}
\partial_r u^*(t,0)=U'(-\tilde c_*t-\xi_0+x_0)<0,\; \partial_r v_*(t,0)=V'(-\tilde c_* t-\xi_0+x_0)>0.
\]
This case is covered by Remark \ref{condition-r=0}, and hence we can still apply the comparison principle to obtain the stated inequalities.
\epf

\begin{rmk}{\rm It is easily seen from the above proof that if instead of \eqref{initial}, we assume that there exists $T>0$ such that $u(T,r)$ and $v(T,r)$ satisfy
the conditions in \eqref{initial}, then the conclusion of Theorem \ref{thm-speed-ub} remains valid.}
\end{rmk}

\end{document}